\documentclass[10pt]{article}

\usepackage{amsmath,amssymb}
\usepackage[a4paper]{geometry}

\newtheorem{Lemme}{Lemma}[section] 
\newtheorem{theorem}[Lemme]{Theorem}
\newtheorem{proposition}[Lemme]{Proposition}
\newtheorem{lemma}[Lemme]{Lemma}   
\newtheorem{corollary}[Lemme]{Corollary}
\newtheorem{remark}[Lemme]{Remark}	
\newtheorem{definition}[Lemme]{Definition}

\def\Box{\leavevmode\vrule height 5pt width 4pt depth 0pt\relax}
\date{\today}

\begin{document}

\begin{center}
{\Large \textbf{Optimal control of steady second grade fluids with a Navier-slip boundary condition}}

\vspace{1 cm}
\textsc{Nadir Arada}{\footnote{Departamento de Matem\'atica, Faculdade de Ci\^encias e Tecnologia da Universidade Nova de Lisboa. E-mail: naar{\char'100}fct.unl.pt.},} 
  \hspace{10mm}
\textsc{Fernanda Cipriano}{\footnote{Departamento de Matem\'atica, Faculdade de Ci\^encias e Tecnologia da Universidade Nova de Lisboa. E-mail: cipriano{\char'100}fct.unl.pt.}}



\end{center}

\vspace{ 8mm}
\noindent \begin{abstract} \noindent In this paper, we first investigate necessary optimality conditions for problems governed
by systems describing the flow of an incompressible second grade fluid. Next, we study the asymptotic behavior of the optimal solution when the viscoelastic parameter  tends to zero, and  prove that the corresponding sequence converges to a solution of an optimal control problem governed by the Navier-Stokes equations whose optimality conditions are recovered by passage to the limit. \vspace{2mm}\vspace{2mm}\vspace{0mm}
\end{abstract} 
\section{Introduction}
\setcounter{equation}{0}
As is well known, the second grade fluid model forms a subclass of 
differential type fluids (also called Rivlin-Ericksen fluids) of complexity 2, and is one of the simplest constitutive models for flows of non-Newtonian fluids that can predict normal stress differences (cf. \cite{RE55} or \cite{NT65}). The Cauchy stress tensor $\mathbb{T}$ for a homogeneous incompressible second grade fluid is given by a constitutive equation of the form
	\begin{equation}\label{cauchy}
	\mathbb{T}=-\pi\mathbb{I}+\nu
	A_{1}(\mathbf y)+
	\alpha_{1}A_2(\mathbf y)
	+\alpha_{2}
	\left(A_{1}(\mathbf y)\right)^{2},\end{equation}
where $\pi$ denotes the hydrodynamic pressure, $\nu$ is the viscosity of the fluid, $\alpha_1$ and $\alpha_2$ are viscoelastic parameters (normal stress moduli), $\mathbf y$ is the velocity field, and $A_1$, $A_2$ are the first two Rivlin-Ericksen tensors defined by 
	$$A_1(\mathbf y)=
	\nabla \mathbf y+\nabla \mathbf y^\top,
	\qquad A_2(\mathbf y)= 
	\tfrac{D{A}_{1}(\mathbf y)}{Dt}+
	{A}_{1}(\mathbf y)\nabla \mathbf y
	+\nabla \mathbf y^\top A_1(\mathbf y)$$
with $\frac{D}{Dt}={\partial }_t+\mathbf y\cdot \nabla$  standing for the material derivative. According to \cite{DF74}, if the fluid modelled by equation (\ref{cauchy}) is to be compatible with thermodynamics  in the sense that all motions of the fluid meet the Clausius-Duhem inequality and the assumption that the specific Helmholtz free energy of the fluid is a minimum in equilibrium, then 
	\begin{equation}\label{thermodynamic}
	\nu\geq 0, \qquad
	\alpha_1\geq 0,\qquad \alpha_1+
	\alpha_2=0.\end{equation}
We refer to \cite{DR95} for a critical and extensive historical review of second-order (and higher order) fluid models and, in particular, for a discussion on the sign of the normal stress moduli. Here we will restraint to the simplified case $\alpha_1+\alpha_2=0$, with $\alpha_1\geq 0$ and $\nu>0$. Setting $\alpha_1=\alpha$ and substituting (\ref{cauchy}) into the balance of linear momentum, we can see that the problem of determining the velocity field $\mathbf{y}$ and the associated pressure $\pi$ satisfying the equations governing the flow of an incompressible second grade fluid reduces to 
	\begin{equation}\label{equation_etat_temps}
	\left\{ \begin{array}{ll}\partial_t\left(\mathbf y-\alpha\Delta  \mathbf y\right)-\nu \Delta  \mathbf y+
	\mathbf{curl}\left( \mathbf y-\alpha\Delta  \mathbf y
	\right)\times  \mathbf y+\nabla \pi= \mathbf u& \mbox{in} \ \Omega,\vspace{2mm} \\
             \mathrm{div} \,  \mathbf y=0& \mbox{in} \ \Omega,\end{array}\right.
	\end{equation}
 where $u$ is the given body force, $\Omega\subset \mathbb{R}^2$ is a bounded domain with boundary $\Gamma$. As this equation is set in dimension two, the vector $\mathbf y$ is written in the form $\mathbf y=(y\equiv (y_1,y_2),0)$ in order to define the curl and the vector product, $\mathbf{curl} \, \mathbf y=(0,0, \mathrm{curl} \,  y)$ with $ \mathrm{curl} \,  y=\tfrac{\partial y_2}{\partial x_1}-\tfrac{\partial y_1}{\partial x_2}$.  Even for this simple but mathematically interesting model, the problem is still difficult since the nonlinear term involves derivatives with higher order than the ones appearing in the viscous term. In the inviscid case ($\nu=0$), the second-grade fluid equations are called $\alpha$-Euler equations. Initially proposed as a regularization of the incompressible Euler equations, they are geometrically significant and have been interpreted as a model of turbulence (cf. \cite{HMR981} and \cite{HMR98}). They also inspired another variant, called the $\alpha$-Navier-Stokes equations that turned out to be very relevant in turbulence modeling (cf. \cite{FHT1}, \cite{FHT2} and the references therein). These equations contain the regularizing term $-\nu\Delta\left(\mathbf y-\alpha \Delta\mathbf y\right)$ instead of $\nu\Delta\mathbf y$, making the dissipation stronger and the problem much easier to solve than in the case of second-grade fluids. When $\alpha=0$, the $\alpha$-Navier-Stokes and the second grade fluid equations are equivalent to the Navier-Stokes equations. \vspace{1mm}\\
Equation (\ref{equation_etat_temps}) can be supplemented with different kinds of boundary conditions. Two of them have been particularly considered in the literature:
	$$\begin{array}{llll}& \mbox{Dirichlet.} &  \quad\mathbf y=0
	 &\quad \mbox{on} \ \Gamma,\vspace{2mm}\\
	&\mbox{Navier-slip.} &\quad \mathbf y\cdot \mathbf n=0,
	\qquad \left(\mathbf n\cdot D\mathbf y\right)\cdot 
	\boldsymbol\tau=0 &\quad 
	\mbox{on} \ \Gamma, \end{array}$$
where $n=(n_1,n_2)$ and $\tau=(-n_2,n_1)$ are the unit normal and tangent vectors, respectively, to the boundary $\Gamma$, and $D\mathbf y=\tfrac{\nabla \mathbf y+\nabla \mathbf y^\top}{2}$ is the symmetric part of the velocity gradient. 
The case of Dirichlet boundary conditions have received a lot of attention. It was systematically studied for the first time in \cite{O81} and \cite{CO84} for both steady and unsteady cases. 
A Galerkin's method in the basis of the eigenfunctions of the operator $\mathbf{curl}(\mathbf{curl}(\mathbf y-\alpha\Delta \mathbf y))$ was especially designed to decompose the problem into a mixed elliptic-hyperbolic type, looking for the velocity $\mathbf y$ as a solution of a Stokes-like system coupled to a transport equation satisfied by $\mathbf{curl}\left( \mathbf y-\alpha\Delta  \mathbf y
	\right)$. Under minimal restrictions on the data, this approach allows the authors to establish the existence of solutions (and automatically recover  $H^3$ regularity) in the steady case, and to prove that the time-dependent version admits
 a unique global solution in the two dimensional case. Much work has been done since these pioneering results and, without ambition for completeness, we cite the extensions in \cite{GS94} and \cite{CG97} where global existence for small initial data in three dimensions was established, the former work using a Schauder fixed point argument while the latter considers the decomposition method on the system of Galerkin equations previously mentionned.\vspace{1mm}\\
The case of second grade fluids with Navier boundary conditions was studied in \cite{BR03}. These conditions are known to deeply modify the properties of the equations, generating additional difficulties related with boundary terms to be correctly handled. In return, some mathematical aspects turned out to be more easily treatable. This is for example the case when studying the controllability of the Navier-Stokes equations (see \cite{C96}). This is also the case when dealing with the inviscid limit of their solutions. Indeed, it is well known that the solution to the Navier-Stokes equations with Navier boundary conditions converges, as $\nu$ tends to zero, to a solution to the Euler equations, while no similar conclusion can be reached when dealing with Dirichlet boundary conditions, responsible for the formation of boundary layers (cf. \cite{CMR98}, \cite{IP06}, \cite{K06}, \cite{LNP05}). Similar considerations apply when analyzing the asymptotic bahavior of the solutions of second-grade fluid equations when the elastic response $\alpha$ and/or the viscosity $\nu$ vanish (cf. \cite{BILN12}, \cite{LT10}, \cite{LNTZ15}).\vspace{1mm}\\
This paper deals with the mathematical analysis of an optimal control problem associated with a steady viscous, incompressible second grade fluid. Control is effected through a distributed mechanical force and the objective is to match the velocity field to a given target field. More precisely, the controls and states are constrained to satisfy the following system of partial differential equations
	\begin{equation}\label{equation_etat}
	\left\{ \begin{array}{ll}-\nu \Delta  \mathbf y+
	\mathbf{curl}\left( \mathbf y-\alpha\Delta  \mathbf y
	\right)\times  \mathbf y+\nabla \pi= \mathbf u& \mbox{in} \ \Omega,\vspace{2mm} \\
             \mathrm{div} \,  \mathbf y=0& \mbox{in} \ \Omega,\vspace{2mm}\\
	\mathbf y\cdot \mathbf n=0,\qquad \left(\mathbf n\cdot D\mathbf 
	y\right)\cdot \boldsymbol\tau=0&
	\mbox{on} \ \Gamma\end{array}\right.\end{equation}
and the optimal control problem reads as 
	$$(P_\alpha) 
	 \left\{\begin{array}{ll}\mbox{minimize} & \displaystyle J(u,y)=
        \tfrac{1}{2}\int_\Omega\left|y-y_d\right|^2\,dx+\tfrac{\lambda}{2}
	\int_\Omega\left|u\right|^2\,dx\vspace{1mm}\\
	\mbox{subject to} & (u,y)\in U_{ad}\times H^3(\Omega) \ \mbox{satisfies} \
 (\ref{equation_etat}) \ \mbox{for some} \ \pi\in L^2(\Omega),\end{array}\right.$$
where  $\lambda\geq 0$, $y_d$ is some desired velocity field and $U_{ad}$, the set of admissible controls, is a nonempty closed convex subset of 
$H(\mathrm{curl})=\left\{v\in L^2(\Omega)\mid \mathrm{curl}\, u\in L^2(\Omega) \right\}$. \vspace{1mm}\\
It is well known when dealing with the optimality conditions for problems governed by highly nonlinear equations, that proving 
the G\^ateaux differentiability of the control-to-state mapping is not an easy task (cf. \cite{A12}, \cite{A13}, \cite{A14}, \cite{S05}, \cite{WR10}).  The main issues are similarly connected with the solvability of the corresponding linearized and adjoint equations and are closely related with the regularity of the
coefficients in the main part of the associated differential operators. As already mentioned, the choice of the special Galerkin basis used to study the state equation  is optimal in the sense that it allows us to prove the existence of regular solutions with minimal assumptions on the data. However, this approach is not appropriate to study the linearized and adjoint equations. Its main drawback lies in the fact that it  {\it automatically imposes} the derivation of a  $H^3$ estimate, which is only guaranteed if the coefficient $\mathrm{curl}\left(\mathbf y-\alpha \Delta \mathbf y\right)$  appearing in these equations belongs to $H^1(\Omega)$, i.e. if the state variable $\mathbf y$ belongs $H^4(\Omega)$. This can only be achieved if we consider more regular controls and impose additional restrictions on their size.  This difficulty can be overcome in the case of Navier boundary conditions by expanding the linearized state and adjoint state in a different basis. Unlike the case of Dirichlet boundary conditions, we are able to derive some corresponding $H^2$ a priori estimates and it turns out that theses estimates are sufficient to carry out our analysis and  derive the optimality conditions for $(P_\alpha)$.\vspace{1mm}\\
In this paper, we are also interested in the asymptotic behavior of solutions of $(P_{\alpha})$, when $\alpha$ tends 
to zero. We will prove in particular that
	\begin{equation}\label{stability}\lim_{\alpha\rightarrow 0^+}\min(P_\alpha)=\min(P_0),\end{equation}
 where $(P_0)$ is the optimal control problem governed by the Navier-Stokes equations and defined by
	$$(P_0) \left\{\begin{array}{ll} \mbox{minimize} & J(u,y)
	\vspace{2mm}\\
	\mbox{subject to}  & (u,y)\in U_{ad}\times H^1(\Omega) \ \mbox{such that}
	\vspace{2mm}\\
&\left\{\begin{array}{ll}-\nu \Delta y+
	y\cdot \nabla y+\nabla \pi=u& \quad 
	\mbox{in} \ \Omega,\vspace{2mm}\\
	\mathrm{div}\, y=0& \quad 
	\mbox{in} \ \Omega,\vspace{2mm}\\
	y\cdot n=0,\qquad \left(n\cdot Dy\right)\cdot \tau=0& \quad 
	\mbox{on} \ \Gamma.\end{array}\right.\end{array}
\right.$$
To obtain such a result, we first establish that the sequence of solutions $(y_\alpha)_\alpha$  of (\ref{equation_etat}) converges to $y$, a solution of the Navier-Stokes equation, when $\alpha$ tends to zero. Next we prove that if $(\bar u_\alpha,\bar y_\alpha)$ is a solution to the problem $(P_\alpha)$ then the sequence $(\bar u_\alpha,\bar y_\alpha)_\alpha$ converges  to a solution $(\bar u_0,\bar y_0)$ of $(P_0)$. Another aspect concerns the necessary optimality conditions. To study the asymptotic behavior of these conditions, we analyze the adjoint equations for $(P_\alpha)$ and prove that the sequence of adjoint solutions converges to the solution of the adjoint equation for $(P_0)$. The optimality conditions for $(P_0)$ are then obtained by passing to the limit in the optimality conditions for $(P_\alpha)$\vspace{1mm}\\
The plan of the present paper is as follows. The main results are stated in Section 2. Notation and preliminary results related with the nonlinear terms are given in Section 3. Section 4 is devoted to the existence and uniqueness results for the state and the linearized state equation and to the derivation of the corresponding estimates. In Section 5, we analyze the Lipschitz continuity and the G\^ateaux differentiability of  the control-to-state mapping and we consider the solvability of the adjoint equation in Section 6. Finally, the proof of the main results are given in Section 7. 

\section{Statement of the main results}
\setcounter{equation}{0}
We first establish the existence of solutions for $(P_\alpha)$ and derive the corresponding necessary optimality conditions. 
\begin{theorem} \label{main_existence}Assume that $U_{ad}$ is bounded in $H(\mathrm{curl})$.
Then problem $(P_\alpha)$ admits at least a solution.
\end{theorem}
\begin{theorem} \label{main_1} Let $(\bar{ u}_\alpha,\bar{y}_\alpha)$ be a solution of $(P_\alpha)$. There exists a positive constant 
$\bar \kappa$ depending only on $\Omega$ such that if the following condition holds
	\begin{equation}\label{control_constraint}
	\bar\kappa \left(\left\|\bar u_\alpha\right\|_2+\alpha
	\left\|\mathrm{curl}\, \bar u_\alpha\right\|_2\right)<\nu^2\end{equation}
then there exists a unique $\bar{p}_\alpha\in H^2(\Omega)$ such that 
	\begin{equation}\label{adj_opt_eq_alpha}
	\left\{ \begin{array}{ll}-\nu \Delta  \bar{\mathbf p}_\alpha-
	\mathbf{curl}\,\sigma\left(\bar{\mathbf y}_\alpha\right)\times
	\bar{\mathbf p}_\alpha+\mathbf{curl}
	\left(\sigma\left(\bar{\mathbf y}_\alpha\times \bar{\mathbf p}_\alpha\right)\right)+\nabla \pi= \bar{\mathbf y}_\alpha-\mathbf y_d& \mbox{in} \ \Omega,\vspace{2mm} \\
             \mathrm{div} \,  \bar{\mathbf p}_\alpha=0& \mbox{in} \ \Omega,\vspace{2mm}\\
                \bar{\mathbf p}_\alpha\cdot \mathbf n=0,\qquad \left(\mathbf n\cdot D\bar{\mathbf p}_\alpha\right)\cdot \boldsymbol\tau=0& \mbox{on} \ \Gamma,\end{array}\right.\end{equation}
and
	\begin{equation}\label{opt_control_alpha}\left(\bar{ p}_\alpha+\lambda\bar{ u}_\alpha,v-\bar{ u}_\alpha\right)\geq 0 \qquad 
	\mbox{for all} \ v\in U_{ad}.\end{equation}
\end{theorem}
Next, we consider the asymptotic analysis of the optimal control $(P_\alpha)$. We first prove that if $(\bar u_\alpha,\bar y_\alpha, \bar p_\alpha)$ satisfies the optimality condition, then a cluster point (for an appropriate topology) is admissible for the problem $(P_0)$ and satisfies the corresponding optimality conditions. In Corollary \ref{corollary_assymp} below, we finally establish the stability property (\ref{stability}). 
\begin{theorem} \label{assympt_1}Let $(\bar u_\alpha,\bar y_\alpha,\bar p_\alpha)$ satisfy the optimality conditions $(\ref{adj_opt_eq_alpha})$-$(\ref{opt_control_alpha})$, and assume that $(\bar u_\alpha)_\alpha$ is bounded in 
$H(\rm{curl})$. Then\vspace{2mm}\\
$i)$ $(\bar u_\alpha,\bar y_\alpha)$ converges $($up to a subsequence when $\alpha$ tends to zero$)$ to a feasible point $(\bar u_0,\bar y_0)$ for problem $(P_0)$ for the $\mbox{weak}$-$H(\rm{curl})\times H^1(\Omega)$  topology.\vspace{2mm}\\
$ii)$ Assume that $\bar u_\alpha$ satisfies $(\ref{control_constraint})$. Then $\bar p_\alpha$ converges $($up to a subsequence when $\alpha$ tends to zero$)$ for the weak topology of $H^1(\Omega)$ to $\bar p_0$, weak solution of the adjoint equation
	\begin{equation}\label{adjoint_limit}
	\left\{\begin{array}{ll}-\nu \Delta \bar p_0-\bar y_0\cdot 
	\nabla \bar p_0+\left(\nabla \bar y_0\right)^\top \bar p_0+\nabla \pi=\bar y_0-y_d& \quad 
	\mbox{in} \ \Omega,\vspace{2mm}\\
	\mathrm{div}\, \bar p_0=0& \quad 
	\mbox{in} \ \Omega,\vspace{2mm}\\
	\bar p_0\cdot n=0,\qquad \left(n\cdot D\bar p_0\right)
	\cdot \tau=0& \quad 
	\mbox{on} \ \Gamma\end{array}\right.
	\end{equation}
and satisfying the optimality condition
	$$\left(\bar p_0+\lambda\bar u_0,v-\bar u_0\right)\geq 0 \qquad 
	\mbox{for all} \ v\in U_{ad}.$$
\end{theorem}
\begin{corollary} \label{corollary_assymp} In addition to the assumptions of Theorem $\ref{assympt_1}$, let us assume that $(\bar u_\alpha,\bar y_\alpha)$ is a solution of problem $(P_\alpha)$. Then the conclusion of Theorem $\ref{assympt_1}$ holds true with a limit point $(\bar u_0,\bar y_0)$ solution of $(P_0)$.
\end{corollary}
\section{Notation, assumptions and preliminary results}
\setcounter{equation}{0}
\subsection{Functional setting}
Throughout the paper $\Omega$ is a bounded, non-axisymmetric, simply connected domain in $\mathbb{R}^2$. The boundary of $\Omega$ is denoted by $\Gamma$ and is sufficiently regular. For $u, v\in \mathbb{R}^{2}$, we define the scalar product by	$u\cdot v=\sum_{i=1}^2 u_iv_i$.
For $\eta, \zeta\in \mathbb{R}^{2\times 2}$, we define the scalar product by	$\eta:\zeta=\sum_{i,j=1}^2 \eta_{ij}\zeta_{ij}$.
 We will also use the following notation
	$$\left(u,v\right)=\displaystyle\int_\Omega u(x)\cdot v(x)\,dx,\qquad 
	\left(\eta,\zeta\right)=\displaystyle\int_\Omega \eta(x):\zeta(x)\,dx.$$
The space of infinitely differentiable functions with compact support in $\Omega$ will be denoted by ${\cal D}(\Omega)$. The standard Sobolev spaces are denoted by $W^{k,p}(\Omega)$ ($k\in \mathbb{N}$ and $1<p<\infty$), and their norms by $\|\cdot\|_{k,p}$. We set $W^{k,2}(\Omega)\equiv H^k(\Omega)$ and $\|\cdot\|_{k,2}\equiv \|\cdot\|_{H^k}$. 
Since many of the quantities occuring in the paper are vector-valued functions, the notation will be abridged for the sake of brevity and we will omit the space dimension in the function space notation. (The meaning should be clear from the context.) \vspace{1mm}\\
In order to eliminate the pressure in the weak formulation of the state equation, we will work in divergence-free spaces and we introduce the following Hilbert spaces:
$$\begin{array}{ll}
H&=\left\{v\in L^2(\Omega)\mid \mathrm{div} \, v=0 \ \text{ in }
\Omega \ \mbox{ and } \ v\cdot n=0 \ \mbox{ on }\Gamma\right\},
	\vspace{2mm}\\
V&=\left\{v\in H^1(\Omega)\mid \mathrm{div} \, v=0 \
 \mbox{ in } \ \Omega\mbox{ and } \ v\cdot n=0 \ \text{ on } \ \Gamma
 \right\},\vspace{2mm}\\
W&=\left\{v\in V\cap H^2(\Omega)\mid \left(n\cdot Dv\right)\cdot \tau
=0 \ \ \mbox{on} \ \Gamma\right\}.
\end{array}$$
In the sequel, we set 
	$$\sigma(y)=y-\alpha \Delta y \qquad \mbox{for} \ 
	y\in H^2(\Omega),$$
and denote by $\mathbb{P}:L^2(\Omega)\longrightarrow H$, the 
Helmholtz projector in $L^2(\Omega)$. It is well know that $\mathbb{P}$ is a linear bounded operator and that is characterized by the equality $\mathbb P y=\tilde y$, where $\tilde y$ is given by the Helmholtz decomposition 
	$$y=\tilde y+\nabla \phi, \qquad \tilde y\in H 
	\quad \mbox{and} \quad \phi\in H^1(\Omega).$$
Let us now collect some classical inequalities. First, we consider the Poincar\'e inequality
	$$\left\|y\right\|_2\leq S_{2} \left\|\nabla y\right\|_2
	\qquad \mbox{for all} \ y\in V$$
and the Sobolev inequality
	$$\left\|y\right\|_4\leq S_{4} \left\|\nabla y\right\|_2
	\qquad \mbox{for all} \ y\in V.$$
Finally, we recall some important facts about the Korn inequality. The commonly used variants read as follows
	\begin{align}\label{Korn_general}
	&\mbox{For every} \ y\in H^1(\Omega): \qquad  \left\|y\right\|_{H^1}\leq C_K\left(\left\|D y\right\|_2+\left\|y\right\|_2\right),\\
	&\mbox{For every} \ y\in H^1_0(\Omega): \qquad\left\|\nabla y\right\|_{2}\leq C_K\left\|D y\right\|_2,\nonumber\end{align}
where $C_K$ is a positive constant depending only on $\Omega$. If $y$ satisfies the tangency boundary condition $y\cdot n=0$ instead of
the homogeneous Dirichlet boundary condition, then
	\begin{equation}\label{Korn_tangent}\left\|\nabla y\right\|_{2}\leq C_K\left\|D y\right\|_2\quad \Longleftrightarrow \quad 
	\Omega \ \mbox{is not axisymmetric}.\end{equation}
  (See \cite{DV02} for more details.) This interesting characterization highlights an additional issue encountered when dealing with steady flows in axisymmetric domains. Indeed, besides the fact that (\ref{Korn_tangent}) is not satisfied, the general inequality (\ref{Korn_general}) is not well adapted since, unlike the unsteady case, the $L^2$ norm of $y$ in our equations cannot be controlled. 
\subsection{Auxiliary results}
The aim of this section is to present several results that will be useful throughout the paper. The first result is fundamental and 
 deals with a boundary identity related with the Navier-slip boundary conditions.  It states in particular that the trace of $\mathrm{curl}\, y$ is a linear function on $y$. (See Proposition 1 in \cite{BR03}.)
\begin{lemma} \label{curl_trace}Let $y\in W$. Then, the following identity holds
	$$\mathrm{curl}\, y\big |_\Gamma=y\cdot g\big|_\Gamma \qquad \mbox{where} \ g=2\tfrac{\partial n}{\partial \tau},$$
with $\tfrac{\partial }{\partial \tau}=n_1\tfrac{\partial }{\partial x_2}-n_2\tfrac{\partial }{\partial x_1}$.\end{lemma}
The next two lemmas will be useful when dealing with a priori estimates for the state, linearized and adjoint state equations. Although similar to the ones of Lemma 5 in \cite{BI06} and Propositions 3 in \cite{BR03} and Lemma 2.1 in \cite{CG97}, the corresponding proofs will be given with a special attention being paid to the dependance on the parameter $\alpha$.
\begin{lemma} \label{sigma_psigma}Let $y\in W\cap H^3(\Omega)$. Then, the following estimates hold
	\begin{equation}\label{sigma_psigma1}\left\|\sigma(y)-\mathbb{P}\sigma(y)\right\|_2\leq 
	c\alpha\left\|\nabla y\right\|_2,\end{equation}
	\begin{equation}\label{sigma_psigma2}\left\|\sigma(y)-\mathbb{P}\sigma(y)\right\|_{H^1}\leq 
	c\alpha\left\|y\right\|_{H^2},\end{equation}
where $c$ is a positive constant only depending on $\Omega$.
\end{lemma}
{\bf Proof.} From the definition of $\mathbb{P}$, there exists $\phi\in H^1(\Omega)$ such that
	$$\sigma(y)-\mathbb{P}\sigma(y)=\nabla \phi$$
and thus 
	$$\Delta \phi=\mathrm{div}\left(\nabla \phi\right)=
	\mathrm{div}\left(\sigma(y)-\mathbb{P}\sigma(y)\right)=0.$$
On the other hand, by taking into account Lemma \ref{curl_trace} we obtain
	$$\begin{array}{ll}
	\tfrac{\partial \phi}{\partial n}\big|_{\Gamma}&=
	n\cdot \nabla \phi\big|_{\Gamma}=
	n\cdot \left(\sigma(y)-\mathbb{P}\sigma(y)\right)\big|_{\Gamma}\vspace{3mm}\\
	&=n\cdot \sigma(y)\big|_{\Gamma}=-\alpha\, n\cdot \Delta y\big|_{\Gamma}\vspace{3mm}\\
	&=-\alpha\, \mathbf n\cdot \Delta \mathbf y\big|_{\Gamma}
	=\alpha\, \mathbf n\cdot \mathbf{curl}\left(\mathbf{curl}\,\mathbf y\right)\big|_{\Gamma}\vspace{3mm}\\
	&=\alpha\, \tfrac{\partial}{\partial \tau}
	\left(y\cdot g\right)\big|_{\Gamma}.\end{array}$$
 Since $\tfrac{\partial}{\partial \tau}
	\left(y\cdot g\right)$ is well defined, the result follows by using standard trace estimates and the regularity theory for elliptic equations with Neuman boundary conditions.$\hfill\Box$ 
\begin{lemma} Let $y\in W\cap H^3(\Omega)$. Then, the following estimates hold
\begin{equation}\label{yh_sigma}\left\|y\right\|_{H^2}\leq 
	\tfrac{c}{\alpha}\left\|\sigma(y)\right\|_2, \end{equation}
	\begin{equation}\label{yh_curl_sigma} \left\|y\right\|_{H^2}\leq 
	\tfrac{c}{\alpha}\left\|\mathrm{curl}\,\sigma(y)\right\|_{2},
	\end{equation}
	\begin{equation}\label{y2_sigma} \left\|y\right\|_{H^3}\leq 
	c\left\|\mathrm{curl}\,\Delta y\right\|_{2},\end{equation}
where $c$ is a positive constant only depending on $\Omega$.\end{lemma}
{\bf Proof.} The proof is split into three steps. \vspace{2mm}\\
 {\it Step 1. Proof of the first estimate.} Let us first recall that for $f\in H^m(\Omega)$, $m\in \mathbb{N}$, the following problem
	$$\left\{\begin{array}{ll}-\alpha\Delta h+h+\nabla \pi=f& \quad \mbox{in} \ \Omega,\vspace{2mm}\\
	\mathrm{div} \, h=0& \quad\mbox{in} \ \Omega,\vspace{2mm}\\
	h\cdot n=0, \qquad (n\cdot Dh)\cdot \tau=0& \quad\mbox{on} \ \Gamma,\vspace{2mm}\\\end{array}\right.$$
admits a unique (up to a constant for $\pi$) solution $(h,\pi)\in H^{m+2}(\Omega)\times H^m(\Omega)$ (see \cite{S73}). Classical arguments show that
	$$\left\|h\right\|_2^2+2\alpha \left\|Dh\right\|_2^2=\left(f,h\right)\leq \left\|f\right\|_2\left\|h\right\|_2$$
yielding 
	\begin{equation}\label{est_g_2_h1}\left\|h\right\|_2\leq \left\|f\right\|_2 \quad \mbox{and} \quad 
	\left\|Dh\right\|_2\leq \tfrac{1}{\sqrt{2\alpha}}\left\|f\right\|_2.\end{equation}
On the other hand, due to the regularity results for the Stokes system, we have
	$$\left\|h\right\|_{H^2}\leq 
	c\left\|\tfrac{f-h}{\alpha}\right\|_2.$$  
Taking into account (\ref{est_g_2_h1}), we deduce that
	$$\left\|h\right\|_{H^2}\leq \tfrac{c}{\alpha}\left(\left\|f\right\|_2+\left\|h\right\|_2\right)\leq \tfrac{c}{\alpha}
	\left\|f\right\|_2$$
and the claimed result follows by setting $f=y-\alpha\,\Delta y$ and $\pi=0$.\vspace{2mm}\\
 {\it Step 2. Proof of the second estimate.} Since $\mathrm{curl}\,\sigma(y)\in L^2(\Omega)$
and $\nabla\cdot \left(\mathrm{curl}\,\sigma(y)\right)=0$, there exists a unique vector-potential $\psi\in H^1(\Omega)$ such that
	$$\left\{\begin{array}{ll}\rm{curl}\,
	 \psi=\mathrm{curl}\,\sigma(y) &\quad \mbox{in} \ \Omega, \vspace{2mm}\\
	\nabla \cdot \psi=0 & \quad\mbox{in} \ \Omega,\vspace{2mm}\\
	\psi\cdot n=0 & \quad\mbox{on} \ \Gamma\end{array}\right.$$
and 
	\begin{equation}\label{sigma_phi}
	\left\|\psi\right\|_{H^1}\leq 
	c\left\|\mathrm{curl}\,\sigma(y)\right\|_2.\end{equation}
It follows that
	$${\rm curl} \left(y-\alpha 
	\Delta y-\psi\right)=0$$
and there exists $\pi\in L^2(\Omega)$ such that
	$$y-\alpha 
	\Delta y-\psi+\nabla \pi=0.$$
Hence $y$ is the solution of the Stokes system
	$$-\Delta y+\nabla\left(\tfrac{\pi}{\alpha}\right)=\tfrac{1}{\alpha}
	\left(\psi-y\right)$$
and satisfies
	\begin{equation}\label{y_phi}\left\|y\right\|_{H^2}
	\leq \tfrac{c}{\alpha}
	\left\|\psi-y\right\|_2.\end{equation}
Since 
	$$\left(\psi,y\right)=
	\left(y-\alpha \Delta y+\nabla \pi,y\right)=
	\left\|y\right\|_2^2+
	2\alpha \left\|Dy\right\|_2^2,$$
we deduce that
	$$\left\|\psi-y\right\|_2^2=
	\left\|\psi\right\|_2^2-
	\left\|y\right\|_2^2-4\alpha
	 \left\|D y\right\|_2^2\leq
	 \left\|\psi\right\|_2^2.$$
Estimate (\ref{yh_curl_sigma}) is then a direct consequence of (\ref{sigma_phi}) and (\ref{y_phi}). \vspace{2mm}\\
 {\it Step 3. Proof of the third estimate.} Since $\mathrm{curl}\,\Delta y\in L^2(\Omega)$ and $\nabla\cdot \left(\mathrm{curl}\,\Delta y\right)=0$, by arguing as in the second step, we can prove that there exists a unique function $\psi\in H^1(\Omega)$ such that 
	$$-\Delta y -\psi+\nabla \pi=0$$
and 
	$$\left\|y\right\|_{H^3}\leq c\left\|\psi\right\|_{H^1}\leq 
	\left\|\mathrm{curl}\,\Delta y\right\|_2$$
which prove the claimed result.$\hfill \Box$
\begin{remark}According to $(\ref{yh_sigma})$ and $(\ref{y2_sigma})$, the standard norms in $H^2(\Omega)$ and $H^3(\Omega)$ are respectively equivalent to the following norms
	$$\left|\cdot \right|_{H^2}=
	\left(\left\|\cdot \right\|_{H^1}^2+\left\|\sigma(\cdot)
	\right\|_2^2\right)^{\frac{1}{2}}\qquad \mbox{and} \qquad
	\left|\cdot \right|_{H^3}=
	\left(\left\|\cdot \right\|_{H^2}^2
	+\left\|\mathrm{curl}\,\sigma(\cdot)
	\right\|_2^2\right)^{\frac{1}{2}}.$$
Moreover, combining $(\ref{yh_sigma})$ and $(\ref{sigma_psigma1})$, we can easily see that the standard norm in $H^2(\Omega)$ is equivalent to the norm
	$$\left\||\cdot \hspace{-0.8mm}|\right\|_{H^2}=
	\left(\left\|\cdot \right\|_{H^1}^2+\left\|\mathbb{P}\sigma(\cdot)
	\right\|_2^2\right)^{\frac{1}{2}}.$$
\end{remark}
The first identity in the next result is standard and relates the nonlinear term in (\ref{equation_etat}), and similar terms appearing in the  linearized and adjoint state equations, to the classical trilinear form used in the Euler and Navier-Stokes equations and defined by
	$$b(\phi,z,y)=\left(\phi\cdot \nabla z,y\right).$$  
The second identity deals with another term only appearing in the adjoint state equation. Unlike the Dirichlet boundary conditions for which the corresponding proofs are straightforward, the Navier-slip boundary conditions are more delicate to handle and proving that the boundary terms, induced by the performed integrations by parts, are vanishing is not an obvious issue.
\begin{lemma} \label{non_lin_curl}
Let $y,z \in  W\cap H^3(\Omega)$ and $\phi\in  V$. Then
	$$\left( \mathbf{curl}\, \sigma(\mathbf y)
	\times  \mathbf z, \boldsymbol\phi\right)=
	b\left(\phi,z, \sigma(y)\right) 
	-b\left(z,\phi,\sigma(y)\right).$$
Let $y, z$ be in $W\cap H^3(\Omega)$ and $\phi$ be in $W$. Then 
	$$\left( \mathbf{curl}\, \sigma\left(\mathbf y\times \mathbf z\right), \boldsymbol\phi\right)
	=b\left(z,y, \sigma(\phi)\right)-b\left(y,z,\sigma(\phi)\right).$$
\end{lemma}
As will be seen in the sequel, the first identity in Lemma \ref{non_lin_curl} enables us to give an adequate variational setting for the state and linearized state equations. Based on the corresponding definitions, we can establish $H^1$ and $H^3$ a priori estimates for the first equation and $H^1$ and $H^2$ a priori estimates for the second. 
Similarly, combining the two identities in Lemma \ref{non_lin_curl}, we can propose a variational formulation for the adjoint equation and establish a $H^1$ estimate of the corresponding solution. However, this formulation cannot be used to establish higher regularity estimates. To overcome this difficulty, we need to use a more suitable, although equivalent, formulation. 
This is the aim of the next lemma.
\begin{lemma}\label{prop_a_adj_2}
Let $y,z \in  W\cap H^3(\Omega)$ and $\phi\in  V$. Then
	$$\begin{array}{ll}
	\left(\mathbf{curl}\left(\sigma\left(\mathbf y\times \mathbf z\right)\right),\boldsymbol\phi\right)
	&=
	b\left(\sigma(z),y,\phi\right)+b\left(y,\phi,\sigma(z)\right)-b\left(\sigma(y),z,\phi\right)+b\left(z,\sigma(y),\phi\right)\vspace{1mm}\\
	&+ b\left(y, z,\phi\right)-b\left(z,y,\phi\right)-2\alpha\displaystyle\sum_{i=1}^2\left(b\left(
	\tfrac{\partial z}{\partial x_i},\tfrac{\partial y}{\partial x_i},\phi\right)
	-b\left(\tfrac{\partial y}{\partial x_i},\tfrac{\partial z}{\partial x_i},\phi\right)\right).
	\end{array}$$
\end{lemma}
We finish this section by stating a lemma that will be used to establish a uniqueness result for the state equation and to derive $H^1$ a priori estimates for the linearized state equation and the adjoint equation. As already observed, 
some specific difficulties related to the boundary terms arise when considering the Navier-slip conditions.
If they are not vanishing, these terms need to be managed and satisfactorily estimated.
\begin{lemma}\label{rm2}
Let $ y \in  W\cap H^3(\Omega)$ and $ z\in  W$. Then
	$$
	\left|\left( \mathbf{curl}\,  \sigma(\mathbf z)	\times  \mathbf y, \mathbf z\right)\right|
	\leq \left( \kappa_1\left\|D y\right\|_2+ \kappa  
	\alpha\left\|y\right\|_{H^3}\right)\left\|D z\right\|_2^2,$$
where $ \kappa_1=S_{4}^2 C_K^3$ and $\kappa$ is a positive constant only depending on $\Omega$. 
\end{lemma}
The proofs of Lemma \ref{non_lin_curl}, Lemma \ref{prop_a_adj_2} and Lemma \ref{rm2} are given in  the appendix.
\section{State equation}
\setcounter{equation}{0}
\subsection{Existence and uniqueness results for the state equation}
 The state equation can be written in a variational form by taking its scalar product with a test function in $ V$.
\begin{definition} Let $ u\in L^2(\Omega)$. A function $ y\in  W\cap H^3(\Omega)$ is a solution of $(\ref{equation_etat})$ if
	\begin{equation}\label{var_form_state}2\nu\left(Dy,D\phi\right)+\left( \mathbf{curl}\, \sigma(\mathbf y)
	\times  \mathbf y, \boldsymbol\phi\right)=\left( u, \phi\right) \qquad \mbox{for all} \  \phi\in  V.\end{equation}
\end{definition}
Due to Lemma \ref{non_lin_curl}, the nonlinear term in this definition can be understood in the following sense
	$$\begin{array}{ll}\left( \mathbf{curl}\, \sigma(\mathbf y)\times  \mathbf y, \boldsymbol\phi\right)&=
	b\left( \phi, y, \sigma(y)\right)-b\left( y,
	 \phi, \sigma(y)\right)
	\vspace{2mm}\\&=b\left( y, y, \phi\right)
	-\alpha\left(b\left( \phi, y,\Delta  y\right)-b\left( y,
	 \phi,\Delta  y\right)\right).\end{array}$$
Equation (\ref{equation_etat}) was first studied by Cioranescu and Ouazar (\cite{O81}, \cite{CO84}) in the case of Dirichlet boundary conditions and simply connected domains. These authors proved existence and uniqueness of solutions by using Galerkin's method in the basis of the eigenfunctions of the operator $\mathrm{curl}\left(\mathrm{curl}\, \sigma(y)\right)$. This method, designed to decompose the problem into a Stokes-like system for the velocity $y$ and a transport equation for $\mathrm{curl}\,\sigma(y)$, allows to establish the  existence of global solutions with  $H^3$ regularity in the  two dimensional case, and uniqueness and local existence in the three dimensional case. It has been  extented by Cioranescu and Girault \cite{CG97} to prove global existence in time in the three dimensional case and by Busuioc and Ratiu \cite{BR03} to study the case of Navier-slip boundary conditions. \vspace{2mm}\\
The following result deals with existence of a solution. It is proved in \cite{BR03} in the unsteady case and can be easily adapted to the steady case. For the convenience of the reader, the corresponding estimates are derived herafter. 
\begin{proposition} \label{existence_state}Let $ u\in H( \mathrm{curl})$. Then problem $(\ref{equation_etat})$ admits at least one solution $ y\in  W\cap H^3(\Omega)$ and this solution satisfies the following estimates
	\begin{equation}\label{state_est1}
	\left\|D y\right\|_{2}\leq \tfrac{ \kappa_2}{\nu}
	\left\| u\right\|_2,\end{equation}
	\begin{equation}\label{state_est2}
	\left\|\mathrm{curl}\, \sigma(y)
	\right\|_{2}\leq \tfrac{\kappa}{\nu}
	\left(\left\| u\right\|_2+\alpha
	\left\| \mathrm{curl} \, u\right\|_2\right),\end{equation}
	\begin{equation}\label{state_est3}
	\left\| y\right\|_{H^3}\leq 
	\tfrac{\kappa}{\alpha\nu}\left(\left\| u\right\|_2+\alpha
	\left\| \mathrm{curl} \, u\right\|_2\right),
	\end{equation}
where $ \kappa_2=\tfrac{S_{2}\,C_K}{2}$ and $\kappa$ is a positive constant depending only on $\Omega$.
\end{proposition}
{\bf Proof.}  Setting $\phi=y$ in (\ref{var_form_state}) and using the Poincar\'e and the Korn inequalities, we obtain
	$$\begin{array}{ll}2\nu\left\|Dy\right\|_2^2&=\left(u,y\right)-
	\left( \mathbf{curl}\, \sigma\mathbf (\mathbf y)\times  \mathbf y, \mathbf y\right)
	=\left(u,y\right)\vspace{2mm}\\
	&\leq
	 \|u\|_2\|y\|_2\leq S_{2}\,C_K\|u\|_2\left\|Dy\right\|_2\end{array}$$
which gives (\ref{state_est1}). On the other hand, by applying the curl to 
$(\ref{equation_etat})$, we obtain
	$$-\nu \Delta \left(\mathrm{curl} \, y\right)+y\cdot \nabla
	 \mathrm{curl} \,\sigma(y)=\mathrm{curl} \, u$$
yielding
	\begin{equation}\label{transport_state}
	\mathrm{curl} \, \sigma(y)+\tfrac{\alpha}{\nu}\,
	y\cdot \nabla\left(\mathrm{curl} \,\sigma(y)\right)
	=\tfrac{\alpha}{\nu}\,
	\mathrm{curl} \, u+\mathrm{curl} \, y.\end{equation}
Multiplying by $\mathrm{curl} \, \sigma(y)$ and integrating, we get
	$$\begin{array}{ll}\left\|\mathrm{curl}\, \sigma(y)\right\|_{2}^2&
	=-\tfrac{\alpha}{\nu}\left(y\cdot \nabla\left(\mathrm{curl} \,\sigma(y)\right),\mathrm{curl}\, \sigma(y)\right)+\left(\tfrac{\alpha}{\nu}\,\mathrm{curl} \, u+\mathrm{curl} \, y,\mathrm{curl}\, 
	\sigma(y)\right)\vspace{2mm}\\
	&=\left(\tfrac{\alpha}{\nu}\,\mathrm{curl} \, u+\mathrm{curl} \, y,\mathrm{curl}\,
	 \sigma(y)\right)	\vspace{2mm}\\
	&\leq\left(\tfrac{\alpha}{\nu}\left\|\mathrm{curl} \, u\right\|_2+\left\|\mathrm{curl} \, y\right\|_2\right)
	\left\|\mathrm{curl}\, \sigma(y)\right\|_2
	\end{array}$$
and thus
	$$\left\|\mathrm{curl}\, \sigma(y)
	\right\|_{2}\leq \left\| \mathrm{curl} \, y\right\|_2
	+\tfrac{\alpha}{\nu}
	\left\| \mathrm{curl} \, u\right\|_2.$$
This estimate together with (\ref{state_est1}) gives (\ref{state_est2}). Finally, by taking into account $(\ref{y2_sigma})$ and 
(\ref{transport_state}), we get
	$$\begin{array}{ll}\left\|y\right\|^2_{H^3}&\leq 
	c^2\left\|\mathrm{curl} \, \Delta y\right\|_2^2
	=\left(\tfrac{c}{\alpha}\right)^2\left\|\mathrm{curl} \, 
	y-\mathrm{curl} \, \sigma(y)\right\|_2^2\vspace{2mm}\\
	&=\left(\tfrac{c}{\alpha}\right)^2\left(
	\left\|\mathrm{curl} \, y\right\|_2^2
	+\left\|\mathrm{curl} \, \sigma(y)\right\|_2^2-2
	\left(\mathrm{curl} \, y,\mathrm{curl} \, \sigma(y)\right)
	\right)\vspace{2mm}\\
	&=\left(\tfrac{c}{\alpha}\right)^2\left(
	\left\|\mathrm{curl} \, y\right\|_2^2
	-\left\|\mathrm{curl} \, \sigma(y)\right\|_2^2
	+\tfrac{2\alpha}{\nu}
	\left(\mathrm{curl} \, u,\mathrm{curl} \, \sigma(y)\right)
	\right)\vspace{2mm}\\
	&\leq \left(\tfrac{c}{\alpha}\right)^2\left(
	\left\|\mathrm{curl} \, y\right\|_2^2
	+\left(\tfrac{\alpha}{\nu}
	\left\|\mathrm{curl} \, u\right\|_2
	\right)^2\right)\end{array}$$
and thus
	$$\left\|y\right\|_{H^3}\leq \tfrac{c}{\alpha}
	\left(\left\|\mathrm{curl} \, y\right\|_2+\tfrac{\alpha}{\nu}
	\left\|\mathrm{curl} \, u\right\|_2\right).$$
Estimate (\ref{state_est3}) is then a direct consequence of (\ref{state_est1})$\hfill\Box$\vspace{2mm}\\
As in the case of Navier-Stokes equations, uniqueness of the solution is guaranteed under a restriction on the data. Aditional regularity of the solution is obtained under the same restriction for more regular data.  
\begin{proposition}\label{uniqueness_state} 
Assume that $ u\in H(\mathrm{curl})$. There exists a constant $ \kappa_\star$ independent of $\nu$ and $\alpha$ such that if $ u$ satisfies 
	$$
	 \kappa_\star\left(\left\| u\right\|_2+\alpha
	\left\| \mathrm{curl} \, u\right\|_2\right)<\nu^2,$$
then equation $(\ref{equation_etat})$ admits a unique solution $y$. 
\end{proposition}
{\bf Proof.} Assume that $ y_1$ and $ y_2$ are two solutions of (\ref{equation_etat}) corresponding to $ u$ and denote by $ y$ the difference $ y_1- y_2$. By setting $\phi=y$ in the weak formulation (\ref{var_form_state}),  we deduce that
	$$2\nu\left\|Dy\right\|_2^2+\left( \mathbf{curl} \,\sigma(\mathbf y_1)
	\times  \mathbf y_1- \mathbf{curl}\, \sigma(\mathbf y_2)
	\times  \mathbf y_2,\mathbf y\right)=0.$$
Observing that
	$$\mathbf{curl} \,\sigma(\mathbf y_1)
	\times  \mathbf y_1- \mathbf{curl}\, \sigma(\mathbf y_2)
	\times  \mathbf y_2=\mathbf{curl} \,\sigma(\mathbf y_1)
	\times  \mathbf y+\mathbf{curl} \,\sigma(\mathbf y)	\times \mathbf  y_2,$$
and taking into account  Lemma \ref{non_lin_curl}, we deduce that
	$$2\nu\left\|Dy\right\|_2^2+\left( \mathbf{curl} \,\sigma(\mathbf y)
	\times  \mathbf y_2,\mathbf y\right)=0.$$
By taking into account  Lemma \ref{rm2}, (\ref{state_est1}) and
 (\ref{state_est3}), it follows that
	$$\begin{array}{ll}\left\|Dy\right\|_2^2
	&\leq \tfrac{1}{2\nu}\left( \kappa_1\left\|Dy_2\right\|_2+
	\kappa\alpha\left\| y_2
	\right\|_{H^3}\right)\left\|D y\right\|_2^2\vspace{2mm}\\
	&\leq \tfrac{1}{2\nu^2}\left( \kappa_1 \kappa_2\left\|u\right\|_2+
	\kappa\left(\left\| u\right\|_2+\alpha
	\left\| \mathrm{curl} \, u\right\|_2\right)\right)
	\left\|D y\right\|_2^2\vspace{2mm}\\
	&\leq \tfrac{ \kappa_\star}{\nu^2}\left(\left\| u\right\|_2+\alpha
	\left\| \mathrm{curl} \, u\right\|_2\right)
	\left\|D y\right\|_2^2
	\end{array}$$
implying that $ y_1= y_2$ if the given condition is satisfied.   $\hfill\Box$
\begin{remark} \label{remark1} Notice that
	$ \kappa_\star >\tfrac{ \kappa_1 \kappa_2}{2}$. This  
implies that if $u$ satisfies the condition stated in the previous proposition, then the corresponding Navier-Stokes equation has a unique weak solution. 
\end{remark}

\subsection{Linearized state equation}
As usual, in order to derive the optimality conditions of an optimal control problem governed by a nonlinear equation, we neeed to study the solvability of the corresponding linearized equation in an adequate setting. Its solution is involved in the definition of the directional derivative of the control-to-state mapping and is related, through a suitable Green formula, to the adjoint state.\vspace{2mm}\\
Let  $u\in H( \mathrm{curl})$, let $ y\in W\cap H^3(\Omega)$ be a corresponding solution of $(\ref{equation_etat})$ and consider the linear equation 
\begin{equation}\label{linearized}\left\{
  \begin{array}{ll}
    -\nu\Delta  \mathbf z+ \mathbf{curl}\,\sigma(\mathbf z)\times \mathbf y+
	 \mathbf{curl}\,\sigma(\mathbf y)\times \mathbf z+\nabla \pi=\mathbf w&\quad\mbox{in} \ \Omega,\vspace{2mm}\\
	 \mathrm{div} \,  \mathbf z=0&\quad\mbox{in} \ \Omega,\vspace{2mm}\\
    \mathbf z\cdot \mathbf n=0, \qquad 
	\left(\mathbf n\cdot D\mathbf z\right)\cdot \boldsymbol\tau=0&\quad\mbox{on}\ \Gamma,
  \end{array}
\right .\end{equation}
where $w\in L^2(\Omega)$.
\begin{definition}  A function $ z\in  W$ is a solution of $(\ref{linearized})$ if
	\begin{equation}\label{var_lin}2\nu\left(D z,D  \phi\right)+
	\left(\mathbf{curl}\, \sigma(\mathbf z)\times  \mathbf y+\mathbf{curl}\, \sigma(\mathbf y)\times  \mathbf z, 
	\boldsymbol\phi\right)=
	\left( w, \phi\right) \qquad \mbox{for all} \ 
	\phi\in  V.\end{equation}
\end{definition}
In analogy to the state equation, by taking into account the first identity in Lemma \ref{non_lin_curl}, we can rewrite the previous variational formulation as follows:
	$$2\nu\left(D z,D  \phi\right)+
	b\left( \phi, y, \sigma(z) \right)-b\left( y, \phi, \sigma(z) \right)+
	b\left( \phi, z, \sigma(y) \right)-b\left( z, \phi, \sigma(y) \right)=
	\left( w, \phi\right)$$
for all $\phi\in V$.\vspace{2mm}\\
As already mentioned, the choice of the special Galerkin basis used to study the state equation (\ref{equation_etat}) is optimal because it allows us to prove the existence of $H^3$ solutions with minimal assumptions on the data. However, this basis does not seem appropriate to study the solvability of both the linearized equation (\ref{linearized}) and the adjoint state equation. Indeed, after deriving the $H^1$ a priori estimates, this technique will naturally {\it imposes} the derivation of a $L^2$ estimate for 
$\mathrm{curl}\, \sigma(z)$  (and thus $H^3$ for $z$). This term should satisfy the transport equation
	$$\mathrm{curl}\, \sigma(z)+\tfrac{\alpha}{\nu} \,y\cdot \nabla 
	\left(\mathrm{curl}\, \sigma(z)\right)
	+\tfrac{\alpha}{\nu}\,z\cdot \nabla 
	\left(\mathrm{curl}\, \sigma(y)\right)=\tfrac{\alpha}{\nu}
	\,\mathrm{curl}\, w+\mathrm{curl}\, z$$
and in order to obtain the desired estimate,  we need to guarantee that the coefficient $\mathrm{curl}\,\sigma(y)$ appearing in the linearized operator belongs to $H^1(\Omega)$. Unfortunately, this can only be achieved if we consider more regular data in the state equation and impose additional restrictions on their size. \vspace{2mm}\\
It turns out that $H^2$ a priori estimates for the linearized state and the adjoint state are sufficient to carry out our analysis. 
Formally, the natural way to obtain such an estimate for $z$  would be to multiply (\ref{linearized}) by $\sigma(z)$
and to integrate. The main difficulty is then to deal with the pressure term
	$$\left(\nabla \pi,\sigma(z)\right)=-\left(\pi,\mathrm{div}\, \sigma(z)\right)+\int_\Gamma \pi n\cdot \sigma(z)=\int_\Gamma \pi n\cdot \sigma(z)$$
 that we do no know how to estimate and that does not vanish, unless  $\sigma(z)$ is tangent to the boundary. To overcome this difficulty,  we multiply (\ref{linearized}) by $\mathbb{P}\sigma(z)$ (instead of  $\sigma(z)$) 
and take advantage of  the nice properties induced by the Navier-slip boundary conditions, stated in Lemma \ref{curl_trace} and Lemma \ref{sigma_psigma}.
\begin{proposition} \label{ex_uniq_lin}Let $ u\in H(\mathrm{curl})$ and let $ y\in  W\cap H^3(\Omega)$ be a corresponding solution of $(\ref{equation_etat})$. There exists a positive constant $\kappa^\star\geq  \kappa_\star$ only depending on $\Omega$, such that if 
	$$\kappa^\star 
	\left(\left\| u\right\|_2+\alpha
	\left\| \mathrm{curl} \, u\right\|_2\right)<\nu^2$$
then equation $(\ref{linearized})$ admits a unique solution $ z\in  W$. Moreover, the following estimates hold
	\begin{equation}\label{est_H1_z}	
	\left(1-\tfrac{\kappa^\star}{\nu^2}\left(\left\| u\right\|_2+\alpha
	\left\| \mathrm{curl} \, u\right\|_2\right)\right)\left\|Dz\right\|_2
	\leq \tfrac{\kappa}{\nu} \left\| w\right\|_2,\end{equation}
	\begin{equation}\label{est_V2_z}\left(1-\tfrac{\kappa^\star}{\nu^2}\left(\left\| u\right\|_2+\alpha
	\left\| \mathrm{curl} \, u\right\|_2\right)\right)\left\| \sigma(z)\right\|_2\leq \kappa\left(\left(1+\alpha
	\right)
	\left\|Dz\right\|_2
	+\tfrac{\alpha}{\nu}\|w\|_2\right),\end{equation} 
where $\kappa$ is a positive constant depending only on $\Omega$.
\end{proposition}
The proof of Proposition \ref{ex_uniq_lin} is split into three steps. We first establish the existence of an approximate solution and a first estimate in $ H^1(\Omega)$. Next, we derive an estimate in $H^2(\Omega)$ that allows the passage to the limit. \vspace{2mm}\\
The solution of $(\ref{linearized})$ is constructed by means of Galerkin's discretization, by expanding the linearized state $z$ in a suitable basis considered by Clopeau {\it et}. {\it al.} to study the Navier-Stokes equations with Navier-slip boundary conditions. Following \cite{CMR98}, there exists a set of eigenfunctions $\left(e_j\right)_j\subset  H^3(\Omega)$ of the problem
	\begin{equation}\label{eigen_functions}
	\left\{\begin{array}{ll}-\Delta e_j+\nabla \pi_j=
	\lambda_je_j &\quad \mbox{in} \ \Omega,\vspace{2mm}\\
	\mathrm{div} \, e_j=0& \quad\mbox{in} \ \Omega,\vspace{2mm}\\
	e_j\cdot n=0, \qquad 
	\left(n\cdot De_j\right)\cdot\tau=0&\quad \mbox{on} \ \Gamma\end{array}\right.
	\end{equation}
with 
	$$0<\lambda_1<\cdots<\lambda_j<\cdots \longrightarrow +\infty.$$
The functions $e_j$ form an orthonormal basis in $H$. The approximate problem is defined by
	\begin{equation}\label{faedo_galerkin_lin}
	\left\{\begin{array}{ll}\mbox{Find} \ 
	 z_m=\displaystyle\sum_{j=1}^m \zeta_{j}
	 e_j \ \mbox{solution}, \ \mbox{for} \ 1\leq j\leq m, \ 
	\mbox{of}\vspace{1mm}\\
	2\nu\left(Dz_m,De_j\right)+
	\left( \mathbf{curl}\,\sigma(\mathbf z_m)\times  \mathbf y+
	 \mathbf{curl}\,\sigma(\mathbf y)\times
	  \mathbf z_m, \mathbf e_j\right)=
	\left(w, e_j\right).\end{array}\right.
	\end{equation}
{\it Step 1. Existence of the discretized solution and a priori $H^1$ estimate.} We prove that the $H^1$ estimate can be derived if $u$ satisfies the condition stated in Proposition \ref{uniqueness_state},  guaranteeing uniqueness of the corresponding state. Let $m$ be fixed and consider $P: \ \mathbb{R}^m \longrightarrow \mathbb{R}^m$ defined by
	$$\left(P \zeta\right)_j=2\nu\left(D z_m,D   e_j\right)+\left(\mathbf{curl}\,\sigma(\mathbf z_m)\times \mathbf y+
	\mathbf{curl}\, \sigma(\mathbf y)\times \mathbf z_m,  \mathbf e_j\right)-
	\left( w,  e_j\right),	$$
where $ z_m=\sum_{j=1}^m\zeta_j   e_j$. The mapping $P$ is obviously continuous. Let us prove that $P(\zeta)\cdot \zeta>0$ if $|\zeta|$ is sufficiently large. Classical arguments together with Lemma \ref{rm2} yield
	\begin{align}\label{P_zeta}{P}(\zeta)\cdot \zeta&=
	2\nu \left\|D  z_m\right\|_2^2+\left(\mathbf{curl}\,\sigma(\mathbf z_m)\times  \mathbf y, \mathbf  z_m\right)-
	\left( w, z_m\right)\\
	&\geq \left(2\nu-\left( \kappa_1\left\|Dy\right\|_2+
	\kappa\alpha\| y\|_{H^3}\right)\right)
	\left\|D  z_m\right\|_2^2
	-\left\| w\right\|_2\left\|z_m\right\|_2\nonumber\\
	&\geq \left(2\nu-\left( \kappa_1\left\|Dy\right\|_2+
	\kappa\alpha\| y\|_{H^3}\right)\right)\tfrac{\lambda_1}{2}
	 \left|\zeta\right|^2-\left\| w\right\|_2
	\left|\zeta\right|\nonumber\\
	&\geq \left(2\nu-\tfrac{2 \kappa_\star}{\nu}
	\left(\left\| u\right\|_2+\alpha
	\left\| \mathrm{curl} \, u\right\|_2\right)
	\right)\lambda_1\left|\zeta\right|^2-
	\left\| w\right\|_2\left|\zeta\right|\nonumber\\
	&\longrightarrow +\infty \qquad \mbox{when} \ 
	|\zeta|\rightarrow +\infty.\nonumber\end{align}
Due to the Brouwer theorem, we deduce that there exists $ \zeta^\ast\in \mathbb{R}^m$ such that ${P}\left( \zeta^\ast\right)=0$ and thus $ z_m=\sum_{j=1}^m\zeta_j^\ast   e_j$ is a solution of problem (\ref{faedo_galerkin_lin}). Due to $(\ref{P_zeta})$ and  Lemma \ref{rm2}, we deduce that
	$$\begin{array}{ll}2\nu \left\|D  z_m\right\|_2^2&=-\left( \mathbf{curl}\,\sigma(\mathbf z_m)\times  \mathbf y,  \mathbf z_m\right)+
	\left( w, z_m\right)\vspace{2mm}\\
	&\leq \left(\left( \kappa_1\left\|D y\right\|_2+
	\kappa\alpha
	\left\| y\right\|_{H^3}\right)
	\left\|D  z_m\right\|_2^2+
	\left\| w\right\|_2
	 \left\|D  z_m\right\|_2\right)\vspace{2mm}\\
	&\leq\tfrac{2\kappa_\ast}{\nu}\left(\left\| u\right\|_2+\alpha
	\left\| \mathrm{curl} \, u\right\|_2\right)
	\left\|D  z_m\right\|_2^2+
	\kappa\left\| w\right\|_2
	 \left\|D z_m\right\|_2
	\end{array}$$
which gives 
\begin{equation}\label{est_H1_zm}\left(
	1-\tfrac{ \kappa_\star}{\nu^2}
	\left(\left\| u\right\|_2+\alpha
	\left\| \mathrm{curl} \, u\right\|_2\right)\right)\left\|Dz_m\right\|_2
	\leq \tfrac{\kappa}{2\nu} \left\| w\right\|_2.\end{equation}
{\it Step 2. A priori $H^2$ estimate.} By taking into account (\ref{eigen_functions}), we have
	$$\mathbb{P}\sigma\left(e_j\right)=\left(1+\alpha\lambda_j\right) e_j$$
and multiplying  (\ref{faedo_galerkin_lin}) by $\left(1+\alpha\lambda_j\right)$, we deduce that
	$$2\nu\left(D z_m,D\left(\mathbb{P}\sigma\left(e_j\right)\right)\right)+
	\left( \mathbf{curl}\, \sigma(\mathbf z_m)\times  \mathbf y+
	\mathbf{curl}\, \sigma(\mathbf y)\times
	 \mathbf z_m, \mathbb{P}\sigma\left(\mathbf e_j\right)\right)=
	\left( w, \mathbb{P}\sigma\left(e_j\right)\right).$$
Recalling that $\mathbb{P}\sigma\left(e_j\right)$ is tangent to the boundary, we have
	$$\begin{array}{ll}2\left(D z_m,D\left(\mathbb{P}\sigma\left(e_j\right)\right)\right)\vspace{2mm}\\
	=\displaystyle-\left(\Delta z_m,\mathbb{P}\sigma\left(e_j\right)\right)-\left(\mathbb{P}\sigma\left(e_j\right),\nabla\left(\mathrm{div} \, z_m\right)\right)+2\int_\Gamma \left(Dz_m\cdot n\right)\cdot \mathbb{P}\sigma\left(e_j\right)
	\vspace{2mm}\\
	=\displaystyle-\left(\Delta z_m,\mathbb{P}\sigma\left(e_j\right)\right)\end{array}$$
and thus
	$$-\nu\left(\Delta z_m,\mathbb{P}\sigma\left(e_j\right)\right)+
	\left( \mathbf{curl}\, \sigma(\mathbf z_m)\times  \mathbf y+
	\mathbf{curl}\, \sigma(\mathbf y)\times
	  \mathbf z_m, \mathbb{P}\sigma\left(\mathbf e_j\right)\right)=
	\left( w, \mathbb{P}\sigma\left(e_j\right)\right).$$
Multiplying by $\zeta^\ast_j$ and summing, it follows that
	$$-\nu\left(\Delta z_m,\mathbb{P} \sigma(z_m)\right)+
	\left( \mathbf{curl}\, \sigma(\mathbf z_m)\times  \mathbf y+
	\mathbf{curl}\, \sigma(\mathbf y)\times
	  \mathbf z_m, \mathbb{P} \sigma(\mathbf z_m)\right)=
	\left( w, \mathbb{P} \sigma(z_m)\right)$$
and consequently, we have
	\begin{align}\label{est_h2_1_1_lin}\left\| \sigma(z_m)\right\|_2^2&=
	\left(z_m, \sigma(z_m)\right)-\alpha\left(\Delta z_m, \sigma(z_m)\right)\nonumber\\
	&=\left(z_m, \sigma(z_m)\right)-\alpha\left(\Delta z_m, \sigma(z_m)-\mathbb{P} \sigma(z_m)\right)-
	\alpha\left(\Delta z_m,\mathbb{P} \sigma(z_m)\right)\nonumber\\
	&=\left(z_m, \sigma(z_m)\right)-\alpha\left(\Delta z_m, \sigma(z_m)-\mathbb{P} \sigma(z_m)\right)\nonumber\\
	&\quad +\tfrac{\alpha}{\nu}\left(\left(w,\mathbb{P} \sigma(z_m)\right)-
	\left(\mathbf{curl}\,\sigma(\mathbf z_m)\times \mathbf y+\mathbf{curl}\, \sigma(\mathbf y)\times \mathbf z_m,\mathbb{P} \sigma(\mathbf z_m)\right)\right). \end{align}
Standard calculation together with (\ref{sigma_psigma1}) and $(\ref{yh_sigma})$ show that
	\begin{align}\label{est_h2_2_1_lin}
	\left|\left(z_m, \sigma(z_m)\right)-
	\alpha\left(\Delta z_m, \sigma(z_m)-\mathbb{P} \sigma(z_m)
	\right)\right|
	&\leq \left\|z_m\right\|_2\left\|\sigma(z_m)\right\|_2+
	c\left\|\sigma(z_m)\right\|_2
	\left\| \sigma(z_m)-\mathbb{P} \sigma(z_m)\right\|_2\nonumber\\
	&\leq c(1+\alpha)\left\|Dz_m\right\|_2
	\left\|\sigma(z_m)\right\|_2.\end{align}
Similarly, by taking into account (\ref{sigma_psigma2}) and $(\ref{yh_sigma})$, we have
	\begin{align}\label{est_h2_3_1_lin} 
	\left|\left( \mathbf{curl}\, \sigma(\mathbf z_m)\times \mathbf y,
	\mathbb{P} \sigma(\mathbf z_m)\right)\right|
	&\leq
	\left|\left( \mathbf{curl}\, \sigma(\mathbf z_m)\times \mathbf y,\sigma(\mathbf z_m)- \mathbb{P}\sigma(\mathbf z_m)\right)
	\right|+\left|\left( \mathbf{curl}\, \sigma(\mathbf z_m)\times \mathbf y, \sigma(\mathbf z_m)\right)
	\right|	\nonumber\\
	&=\left|b\left(\sigma(z_m)- \mathbb{P}\sigma(z_m),y, \sigma(z_m)\right)
	-b\left(y,\sigma(z_m)- \mathbb{P}\sigma(z_m),\sigma(z_m)\right)\right|\nonumber\\
	&\quad +
	\left|b\left( \sigma(z_m),y, \sigma(z_m)\right)\right|\nonumber\\
	&\leq \left\|y\right\|_{1,\infty}
	\left\|\sigma(z_m)- \mathbb{P}\sigma(z_m)\right\|_{H^1}
	\left\| \sigma(z_m)\right\|_2+\left\|\nabla y\right\|_\infty
	\left\| \sigma(z_m)\right\|_2^2\nonumber\\
	&\leq c\alpha\left\|y\right\|_{1,\infty}
	\left\|z_m\right\|_{H^2}
	\left\| \sigma(z_m)\right\|_2+\left\|\nabla y\right\|_\infty
	\left\| \sigma(z_m)\right\|_2^2\nonumber\\
	&\leq c\left\|y\right\|_{1,\infty}
	\left\| \sigma(z_m)\right\|_2^2
	\end{align}
and
	\begin{align}\label{est_h2_4_1_lin} 
	\left|\left( \mathbf{curl}\, \sigma(\mathbf y)\times \mathbf z_m,
	\mathbb{P} \sigma(\mathbf z_m)\right)\right|&\leq 
	\left\|\mathrm{curl}\, \sigma(y)\right\|_2
	\left\|z_m\right\|_\infty
	\left\|\mathbb{P}\sigma(z_m)\right\|_2\nonumber\\
	& \leq c
	\left\| \mathrm{curl}\, \sigma(y)\right\|_2
	\left\|z_m\right\|_{H^2}
	\left\|\sigma(z_m)\right\|_2\nonumber\\
	& \leq \tfrac{c}{\alpha}
	\left\| \mathrm{curl}\, \sigma(y)\right\|_2
	\left\|\sigma(z_m)\right\|_2^2.
	\end{align}
Combining (\ref{est_h2_1_1_lin})-(\ref{est_h2_4_1_lin}), and taking into account (\ref{state_est2}) and 
	 (\ref{state_est3}), we obtain
	$$\begin{array}{ll}
	\left\| \sigma(z_m)\right\|_2^2&\leq 
	c\left((1+\alpha)\left\|Dz_m\right\|_2
	+\tfrac{\alpha}{\nu}\left\|w\right\|_2\right) 
	\left\|\sigma(z_m)\right\|_2\vspace{2mm}\\
	&+\tfrac{c}{\nu}\left(\alpha \left\|y\right\|_{1,\infty}+\left\| \mathrm{curl}\, \sigma(y)\right\|_2\right)
	\left\| \sigma(z_m)\right\|_2^2\vspace{2mm}\\
	&\leq 
	\kappa\left((1+\alpha)\left\|Dz_m\right\|_2
	+\tfrac{\alpha}{\nu}\left\|w\right\|_2\right) 
	\left\|\sigma(z_m)\right\|_2+\tfrac{ \kappa_3}{\nu^2}\left(\left\| u\right\|_2+\alpha
	\left\| \mathrm{curl} \, u\right\|_2\right)\left\| \sigma(z_m)\right\|_2^2 \end{array}$$
which gives
	\begin{equation}\label{est_V2_zm}
	\left(1-\tfrac{ \kappa_3}{\nu^2}\left(\left\| u\right\|_2+\alpha
	\left\| \mathrm{curl} \, u\right\|_2\right)\right)\left\| \sigma(z_m)\right\|_2\leq \kappa\left(\left(1+\alpha
	\right)
	\left\|Dz_m\right\|_2
	+\tfrac{\alpha}{\nu}\|w\|_2\right).\end{equation}
{\it Step 3. Passing to the limit.} It remains to pass to the limit with respect to $m$. From estimates (\ref{est_H1_zm}) and (\ref{est_V2_zm}), it follows that if the following condition
	$$\kappa^\star 
	\left(\left\| u\right\|_2+\alpha
	\left\| \mathrm{curl} \, u\right\|_2\right)<\nu^2 \qquad \mbox{with} \ \kappa^\star=\max\left( \kappa_\star, \kappa_3\right)$$
 is fulfilled, then there exists a subsequence, still indexed by $m$, and a function $ z\in  W$ such that
	$$ z_m \longrightarrow  z \qquad \mbox{weakly in} \  W.$$
By passing to the limit in (\ref{faedo_galerkin_lin}), we obtain for every $j\geq 1$
	$$2\nu\left(D z,D   e_j\right)+
	b\left(  e_j, y, \sigma(z)\right)-b\left( y,  e_j,
	\sigma(z)\right)+b\left(  e_j, z, \sigma(y) \right)
	-b\left( z,  e_j, \sigma(y)\right)=
	\left( w,  e_j\right)	$$
and by density we prove that $ z$ satisfies the variational formulation. Moreover, $z$ satisfies estimates (\ref{est_H1_z}) and (\ref{est_V2_z}). Finally, since (\ref{linearized}) is linear, the uniquess result is a direct consequence of estimate (\ref{est_H1_z}).$\hfill\Box$
\section{Analysis of the control-to-state mapping}
\setcounter{equation}{0}
\subsection{Lipschitz continuity}
In this section, we derive some useful estimates related with the local Lipschitz continuity of the state with respect to the control variable. More precisely, if $u_1$, $u_2$ are two admissible controls and if $y_1$,  $y_2$ are two corresponding states then we are interested in estimating the difference $y_1-y_2$ with respect to $u_1-u_2$ in adequate topologies. The arguments are similar to those used to study the solvability of the linearized state equation but must take into account that a direct adaptation will lead to the imposition of restrictions on both functions $u_1$ and $u_2$. This is particularly the case when dealing with the $H^2$ estimate for $y_1-y_2$ and would result, when deriving the optimality conditions, in all the admissible controls being restrained. This difficulty is overcome by adapting these arguments in order to restrain only one variable.
\begin{proposition}
Let $ u_1,  u_2\in H(\mathrm{curl})$, and let 
$ y_1,  y_2 \in  W\cap H^3(\Omega)$ be corresponding solutions of $(\ref{equation_etat})$. Then the following estimates hold
	\begin{equation}\label{lipschitz_est_H1}
	\left(1-\tfrac{ \kappa_\star}{\nu^2}\left(\left\| u_2\right\|_2+\alpha
	\left\| \mathrm{curl} \, u_2\right\|_2\right)\right)
	\left\|D(y_1- y_2)
	\right\|_2\leq \tfrac{\kappa}{\nu} \left\| u_1- u_2\right\|_2,\end{equation}
	$$\left(1-\tfrac{ \kappa_4}{\nu^2}\left(\left\| u_2\right\|_2+\alpha
	\left\| \mathrm{curl} \, u_2\right\|_2\right)\right)\left\| \sigma(y_1-y_2)\right\|_2^2$$
	\begin{equation}\label{lipschitz_est_V2}
	\leq \kappa\left(\left(1+\alpha
	+\tfrac{\alpha}{\nu^3} \left(\left\| u_1\right\|_2+\alpha
	\left\| \mathrm{curl} \, u_1\right\|_2\right)^\frac{3}{2}\right)
	\left\|D(y_1-y_2)\right\|_2
	+\tfrac{\alpha}{\nu}\|u_1-u_2\|_2\right)^2,
	\end{equation}
where $\kappa$ and $ \kappa_4$ are positive constants depending only on $\Omega$.
\end{proposition}
{\bf Proof.} The proof is split into two steps.\vspace{2mm}\\
{\it Step 1. A priori $H^1$ estimate.}
It is easy to see that $y=y_1-y_2$ satisfies
\begin{equation}\label{y1-y2}\left\{
  \begin{array}{ll}
    -\nu\Delta  \mathbf y+ \mathbf{curl}\, \sigma(\mathbf y)\times  \mathbf y_2+
	\mathbf{curl}\, \sigma(\mathbf y_1)\times  \mathbf y+\nabla \pi= \mathbf u&\quad\mbox{in} \ \Omega,\vspace{2mm}\\
	 \mathrm{div} \,  \mathbf y=0&\quad\mbox{in} \ \Omega,\vspace{2mm}\\
   \mathbf y\cdot \mathbf n=0, \qquad 
	\left(\mathbf n\cdot D\mathbf y\right)\cdot \boldsymbol\tau=0&\quad\mbox{on}\ \Gamma,
  \end{array}
\right.\end{equation}
where $u=u_1-u_2$. By setting $\phi=y$ in the corresponding weak formulation, we obtain
	$$2\nu\left\|Dy\right\|_2^2+
	 \left( \mathbf{curl}\, \sigma(\mathbf y)\times \mathbf  y_2,  \mathbf y\right)=\left( u, y\right).$$
Due to Lemma \ref{rm2}, (\ref{state_est1}) and
 (\ref{state_est3}), it follows that
	$$2\nu\left\|Dy\right\|_2^2
	\leq \left\| u\|_2\| y\right\|_2+
\left(\kappa_1\left\|Dy_2\right\|_2+
	\kappa\alpha\left\|y_2\right\|_{H^3}\right)
	\left\|D y\right\|_2^2$$
 and (\ref{lipschitz_est_H1}) holds. \vspace{2mm}\\
{\it Step 2. A priori $H^2$ estimate.} Multiplying equation $(\ref{y1-y2})_1$ by $\mathbb{P}\sigma(y)$, we obtain
$$-\nu\left(\Delta y,\mathbb{P} \sigma(y)\right)+
	\left( \mathbf{curl}\, \sigma(\mathbf y)\times  \mathbf y_2+
	\mathbf{curl}\, \sigma(\mathbf y_1)\times
	  \mathbf y, \mathbb{P} \sigma(\mathbf y)\right)=
	\left(u, \mathbb{P} \sigma(y)\right).$$
Arguing as in the proof of Proposition \ref{ex_uniq_lin}, we can prove that
	\begin{align}\label{est_h2_1_1_lip}\left\| \sigma(y)\right\|_2^2
	&=\left(y, \sigma(y)\right)-\alpha\left(\Delta y, \sigma(y)-\mathbb{P} \sigma(y)\right)\nonumber\\
	&\quad +\tfrac{\alpha}{\nu}\left(\left(u,\mathbb{P} \sigma(y)\right)-
	\left(\mathbf{curl}\,\sigma(\mathbf y)\times \mathbf y_2+\mathbf{curl}\, \sigma(\mathbf y_1)\times \mathbf y,\mathbb{P} \sigma(\mathbf y)\right)\right) \end{align}
with
	\begin{equation}\label{est_h2_2_1}
	\left|\left(y, \sigma(y)\right)-
	\alpha\left(\Delta y, \sigma(y)-\mathbb{P} \sigma(y)
	\right)\right|\leq c(1+\alpha)\left\|Dy\right\|_2\left\|\sigma(y)\right\|_2\end{equation}
and 
	\begin{equation}\label{est_h2_3_1_lip} 
	\left|\left( \mathbf{curl}\, \sigma(\mathbf y)\times \mathbf y_2,
	\mathbb{P} \sigma(\mathbf y)\right)\right|\leq c\left\|y_2\right\|_{1,\infty}
	\left\| \sigma(y)\right\|_2^2.\end{equation}
Finally, the interpolation inequalities
	$$\left\|y\right\|_{4}\leq c\left\|y\right\|_2^\frac{1}{2}\left\|y\right\|_{H^1}^{\frac{1}{2}}, \qquad 
	\left\|y\right\|_\infty\leq c 
	\left\| y\right\|_2^{\frac{1}{3}} 
	\left\| y\right\|_{1,4}^{\frac{2}{3}}$$
together with (\ref{yh_sigma}) yield
	\begin{align}\label{est_h2_4_1_lip} 
	\left|\left( \mathbf{curl}\, \sigma(\mathbf y_1)\times \mathbf y,
	\mathbb{P} \sigma(\mathbf y)\right)\right|&\leq 
	\left\|\mathrm{curl}\, \sigma(y_1)\right\|_2
	\left\|y\right\|_\infty
	\left\|\mathbb{P}\sigma(y)\right\|_2\nonumber\\
	& \leq c
	\left\| \mathrm{curl}\, \sigma(y_1)\right\|_2
	 \left\|y\right\|_{H^1}^{\frac{2}{3}}
	\left\|y\right\|_{H^2}^{\frac{1}{3}}
	\left\|\sigma(y)\right\|_2\nonumber\\
	& \leq c
	\left\| \mathrm{curl}\, \sigma(y_1)\right\|_2
	 \left\|Dy\right\|_2^{\frac{2}{3}}
	\left\|y\right\|_{H^2}^{\frac{1}{3}}
	\left\|\sigma(y)\right\|_2\nonumber\\
	& \leq \tfrac{c}{\alpha^{\frac{1}{3}}}
	\left\| \mathrm{curl}\, \sigma(y_1)\right\|_2
	 \left\|Dy\right\|_2^{\frac{2}{3}}
	\left\|\sigma(y)\right\|_2^{\frac{4}{3}}.
	\end{align}
Combining (\ref{est_h2_1_1_lip})-(\ref{est_h2_4_1_lip}), we obtain
$$\begin{array}{ll}
	\left\| \sigma(y)\right\|_2^2&\leq 
	c\left((1+\alpha)\left\|Dy\right\|_2
	+\tfrac{\alpha}{\nu}\left\|u\right\|_2\right) 
	\left\|\sigma(y)\right\|_2+
	\tfrac{c\alpha}{\nu}\left\|y_2\right\|_{1,\infty}
	\left\| \sigma(y)\right\|_2^2\vspace{2mm}\\
	& \ \ 
	+\tfrac{c\alpha^{\frac{2}{3}}}{\nu}\left\| \mathrm{curl}\, \sigma(y_1)\right\|_2
	 \left\|Dy\right\|_2^{\frac{2}{3}} \left\|\sigma(y)
	\right\|_2^{\frac{4}{3}}\end{array}$$
and by  using the Young inequality, it follows that
	$$\left\| \sigma(y)\right\|_2^2\leq c\left(\left((1+\alpha)^2
	+\tfrac{\alpha^2}{\nu^3}
	\left\| \mathrm{curl}\,\sigma(y_1)\right\|_2^3\right)
	\left\|Dy\right\|_2^2
	+\left(\tfrac{\alpha}{\nu}\|u\|_2\right)^2
	+\tfrac{\alpha}{\nu}\left\|y_2\right\|_{1,\infty}
	\left\| \sigma(y)\right\|_2^2\right).$$
Hence, by taking into account the estimates (\ref{state_est1})-(\ref{state_est3}), we deduce that 
	$$\begin{array}{ll}
	\left(1-\tfrac{ \kappa_4}{\nu^2}\left(\left\| u_2\right\|_2+\alpha
	\left\| \mathrm{curl} \, u_2\right\|_2\right)\right)\left\| \sigma(y)\right\|_2^2\vspace{2mm}\\
	\leq \kappa\left(\left(1+\alpha
	+\tfrac{\alpha}{\nu^3} \left(\left\| u_1\right\|_2+\alpha
	\left\| \mathrm{curl} \, u_1\right\|_2\right)^\frac{3}{2}\right)
	\left\|Dy\right\|_2
	+\tfrac{\alpha}{\nu}\|u\|_2\right)^2\end{array}$$
which gives (\ref{lipschitz_est_V2}).$\hfill\Box$
\subsection{G\^ateaux differentiability}
At this stage, we are able to study the G\^ateaux-differentiability of the control-to-state mapping. 
\begin{proposition}\label{taylor}
Let $\rho$ be such that $0<\rho<1$, and
 $ u, w$ be in $H(\mathrm{curl})$. Set  $ u_\rho= u+\rho  w$, and let $y$ and $ y_{\rho}$ be solutions of 
$(\ref{equation_etat})$ corresponding to $u$ and $u_\rho$, respectively. There exists a positive constant $\kappa^{\star\star}\geq  \kappa_\star$ only depending on $\Omega$, such that if 
	$$\kappa^{\star\star}
	\left(\left\| u\right\|_2+\alpha
	\left\| \mathrm{curl} \, u\right\|_2\right)<\nu^2$$
then we have
	$$ y_{\rho}= y+\rho 
	 z+\rho r_\rho \qquad 
	\mbox{with} \ \lim_{\rho\rightarrow 0}
	\left\|D r_\rho\right\|_2=0,$$
and
	$$J\left(u_\rho,y_\rho\right)=J\left(u,y\right)+\rho 
	\left(\left(z,y-y_d\right)+\lambda(u,w)\right)
	+o(\rho),$$
where $ z\in  W$ is a solution of $(\ref{linearized})$ corresponding to $(u,y)$.
\end{proposition}
{\bf Proof.} Easy calculation shows that $z_{\rho}=\frac{ y_{\rho}- y}{\rho}$ satisfies
	$$-\nu \Delta  \mathbf z_\rho+ \mathbf{curl}\, 
	\sigma\left(\mathbf z_\rho\right)\times \mathbf y+ \mathbf{curl}\,\sigma\left(\mathbf y_\rho\right)\times \mathbf z_\rho+\nabla\pi_\rho= \mathbf w.$$ 
Let $ z\in  W$ be the solution of (\ref{linearized}).	
Then $r_\rho= z_\rho- z$ satisfies
	$$-\nu \Delta \mathbf r_\rho+ 
	\mathbf{curl}\,\sigma\left(\mathbf r_\rho\right)\times \mathbf y
	+ \mathbf{curl}\,\sigma\left(\mathbf y_\rho\right)\times
	 \mathbf r_\rho
	+ \mathbf{curl}\,\sigma\left(\mathbf y_\rho-\mathbf y\right)
	\times \mathbf z
	+\nabla\left(\pi_\rho-\pi\right)=0.$$
Multiplying this equation by $r_\rho$, we obtain	
	\begin{equation}\label{energy_r}2\nu \left\|Dr_\rho\right\|_2^2+
	\left(\mathbf{curl}\,\sigma\left(\mathbf r_\rho\right)\times  \mathbf y+ \mathbf{curl}\,\sigma\left(\mathbf y_\rho\right)
	\times \mathbf r_\rho+\mathbf{curl}\,\sigma\left(\mathbf y_\rho-\mathbf y\right)
	\times \mathbf z,\mathbf r_\rho\right)=0.\end{equation}
It is easy to verify that
	$$\left( \mathbf{curl}\,\sigma\left(\mathbf y_\rho\right)
	\times \mathbf r_\rho,\mathbf r_\rho\right)=
	b\left(r_\rho,r_\rho,\sigma\left(y_\rho\right)\right)-
	b\left(r_\rho,r_\rho,\sigma\left(y_\rho\right)\right)=0.$$
Moreover, by taking into account Lemma \ref{rm2} and estimates (\ref{state_est1})-(\ref{state_est3}), we get
	\begin{equation}\left|\left(\mathbf{curl}\,
	\sigma\left(\mathbf r_\rho\right)\times  \mathbf y,\mathbf 
	r_\rho\right)\right|
	\leq \tfrac{2 \kappa_\star}{\nu}\left(\left\| u\right\|_2+\alpha
	\left\| \mathrm{curl} \, u\right\|_2\right)
	\left\|Dr_\rho\right\|_2^2
	\end{equation}
and similarly
	\begin{align}\label{energy_r4}\left|\left(\mathbf{curl}\,\sigma\left(\mathbf y_\rho-\mathbf y\right)
	\times \mathbf z,\mathbf r_\rho\right)\right|&=
	\left|b\left(r_\rho,z,\sigma\left(y_\rho-y\right)\right)-
	b\left(z,r_\rho,\sigma\left(y_\rho-y\right)\right)\right|
	\nonumber\\
	&\leq \left(\left\|r_\rho\right\|_4 \left\|\nabla z\right\|_4
+\left\|z\right\|_\infty \left\|\nabla r_\rho\right\|_2\right)
	\left\|\sigma\left(y_\rho-y\right)\right\|_2\nonumber\\
	&\leq c\left\|z\right\|_{H^1} \left\|Dr_\rho\right\|_2
	\left\|\sigma\left(y_\rho-y\right)\right\|_2
	\end{align}
Combining (\ref{energy_r})-(\ref{energy_r4}), we deduce that
	\begin{equation}\label{nr_1}
	\left(1-\tfrac{ \kappa_\star}{\nu^2}\left(\left\| u\right\|_2+\alpha
	\left\| \mathrm{curl} \, u\right\|_2\right)\right)
	\left\|Dr_\rho\right\|_2\leq \tfrac{\kappa}{\nu}
	\left\|z\right\|_{H^1} 
	\left\|\sigma\left(y_\rho-y\right)\right\|_2.\end{equation}
On the other hand, due to (\ref{lipschitz_est_H1})-(\ref{lipschitz_est_V2}), we have
	\begin{equation}\label{nr_2}\left(1-\tfrac{ \kappa_\star}{\nu^2}\left(\left\| u\right\|_2+\alpha
	\left\| \mathrm{curl} \, u\right\|_2\right)
	\right)\left\|D(y_\rho- y)\right\|_2\leq
	 \tfrac{\kappa}{\nu} \rho \left\|w\right\|_2,\end{equation}
and
	\begin{align}\label{nr_3}&\left(1-\tfrac{ \kappa_4}{\nu^2} \left(\left\| u\right\|_2+\alpha
	\left\| \mathrm{curl} \, u\right\|_2\right)
	\right)\left\| \sigma(y_\rho-y)\right\|_2^2\nonumber\\
	&\leq 	\kappa\left(\left(1+\alpha
	+\tfrac{\alpha}{\nu^3} \left(\left\| u_\rho\right\|_2+\alpha
	\left\| \mathrm{curl} \, u_\rho\right\|_2
	\right)^\frac{3}{2}\right)\left\|D(y_\rho-y)\right\|_2
	+\tfrac{\alpha}{\nu}\rho\|w\|_2\right)^2.\end{align}
From (\ref{nr_1})-(\ref{nr_3}), we deduce that if 
	$$\kappa^{\star\star}
	\left(\left\| u\right\|_2+\alpha
	\left\| \mathrm{curl} \, u\right\|_2\right)<\nu^2 \qquad \mbox{with} \ 
	\kappa^{\star\star}=\max\left( \kappa_\ast, \kappa_4\right)$$
then	
	$$\lim_{\rho\rightarrow 0}\left\|D r_\rho\right\|_2=0.$$
The second assertion can be easily proved. $\hfill\Box$
\section{Adjoint equation}
\setcounter{equation}{0}
Let  $u\in H( \mathrm{curl})$ and let $ y\in W\cap H^3(\Omega)$ be a corresponding solution of (\ref{equation_etat}). The aim of this section is to study the solvability of the adjoint state equation defined by
\begin{equation}\label{adjoint}\left\{
  \begin{array}{ll}
    -\nu\Delta  \mathbf p-  \mathbf{curl}\,\sigma(\mathbf y)\times  \mathbf p +\mathbf{curl}\left(\sigma\left(\mathbf y\times \mathbf p\right)\right)+\nabla \pi=\mathbf f&\quad\mbox{in} \ \Omega,\vspace{2mm}\\
	 \mathrm{div} \,  \mathbf p=0&\quad\mbox{in} \ \Omega,\vspace{2mm}\\
    \mathbf p\cdot \mathbf n=0, \qquad 
	\left(\mathbf n\cdot D\mathbf p\right)\cdot \boldsymbol\tau=0&\quad\mbox{on}\ \Gamma,
  \end{array}
\right .\end{equation}
where $f\in L^2(\Omega)$.
The two identities in Lemma \ref{non_lin_curl} motivates the following variational formulation.
\begin{definition} \label{var_adj_2}A function $p\in  W$ is a solution of $(\ref{adjoint})$ if
	\begin{equation}\label{form_var_lin_adj}
	2\nu\left(D p,D  \phi\right)+
	b\left( p,\phi, \sigma(y) \right)-b\left( \phi,p, \sigma(y) \right)+
	b\left(p, y, \sigma(\phi) \right)-b\left( y,p, \sigma(\phi) \right)=
	\left(f, \phi\right)\end{equation}
for all $\phi\in W$.
\end{definition}
This formulation allows us to relate the adjoint state $p$ to the solution $z$ of the linearized equation and is particularly suited to derive the necessary optimality conditions. Indeed, if $p$ is a solution of $(\ref{adjoint})$ in the sense of the previous definition and if $z$ is a solution of (\ref{linearized}), then $z$ is an admissible test function for (\ref{form_var_lin_adj}), $p$ is an admissible test function for (\ref{var_lin}) and
	$$\left(f,z\right)=\left(w,p\right).$$
Let us now observe that the previous reasoning is also valid if $p$ is less regular and belongs to $V$. The optimality conditions can be derived similarly, but the uniqueness of the solution for (\ref{adjoint}) cannot be guaranteed and the corresponding estimate in $H^2(\Omega)$ cannot be established. \vspace{1mm}\\
These considerations can be addressed when considering the approximated problem. 
The solution of $(\ref{adjoint})$ is constructed using the Galerkin's method defined in the previous section. Existence of an approximate solution and a corresponding a priori $H^1$ estimate can first be established by taking into account the formulation stated in Definition \ref{var_adj_2}. This estimate is sufficient to pass to the limit and prove the existence of $p\in V$ satisfying $(\ref{form_var_lin_adj})$. In order to guarantee the uniqueness of this solution, we need to prove that it belongs to $W$ and, thus, to establish first a $H^2$ a priori estimate for the approximate solution. Unlike the linearized equation, the same variational formulation cannot be used to establish both estimates. This issue is related with the term  $\mathbf{curl}\left(\sigma\left(\mathbf y\times \cdot\right)\right)$ and can be  overcome by considering an equivalent formulation, motivated by Lemma \ref{prop_a_adj_2} and more consistent with our objective.  
 \begin{proposition} \label{ex_uniq_adj}Let $ u\in H(\mathrm{curl})$ and let $ y\in  W\cap H^3(\Omega)$ be a corresponding solution of $(\ref{equation_etat})$. There exists a positive constant $ \kappa_{\star\star}\geq  \kappa_\star$ only depending on $\Omega$, such that if 
	$$ \kappa_{\star\star}
	\left(\left\| u\right\|_2+\alpha
	\left\| \mathrm{curl} \, u\right\|_2\right)<\nu^2$$
then equation $(\ref{adjoint})$ admits a unique solution $p\in  W$. Moreover, the following etimates hold
	\begin{equation}\label{est_H1_p}	
	\left(1-\tfrac{ \kappa_{\star\star}}{\nu^2}\left(\left\| u\right\|_2+\alpha
	\left\| \mathrm{curl} \, u\right\|_2\right)\right)\left\|Dp\right\|_2
	\leq \tfrac{\kappa}{\nu} \left\|f\right\|_2,\end{equation}
	\begin{equation}\label{est_V2_p}\left(1-\tfrac{ \kappa_{\star\star}}{\nu^2}\left(\left\| u\right\|_2+\alpha
	\left\| \mathrm{curl} \, u\right\|_2\right)\right)\left\| \sigma(p)\right\|_2\leq \kappa\left(\left(1+\alpha\right)
	\left\|Dp\right\|_2
	+\tfrac{\alpha}{\nu}\|f\|_2\right),\end{equation}
where $\kappa$ is a positive constant only depending on $\Omega$.
\end{proposition}
The approximate problem is defined by
	\begin{equation}\label{faedo_galerkin_adj}
	\left\{\begin{array}{ll}\mbox{Find} \ 
	 p_m=\displaystyle\sum_{j=1}^m \zeta_{j}
	 e_j \ \mbox{solution}, \ \mbox{for} \ 1\leq j\leq m, \ 
	\mbox{of}\vspace{3mm}\\
	2\nu\left(Dp_m,De_j\right)-\left( \mathbf{curl}\,\sigma(\mathbf y)\times  \mathbf p_m,e_j\right)+\left(\mathbf{curl}\left(\sigma\left(\mathbf y\times \mathbf p_m\right)\right),e_j\right)=
	\left(f, e_j\right),\end{array}\right.
	\end{equation}
where $(e_j)_j\subset H^3(\Omega)$ is the set of the eigenfunctions, solutions of (\ref{eigen_functions}). The proof of Proposition \ref{ex_uniq_adj} is split into three steps. We first establish existence of an approximate solution and a corresponding apriori estimate in $H^1(\Omega)$ , we next derive a corresponding estimate in $H^2(\Omega)$  and we finally pass to the limit.\vspace{2mm}\\
{\it Step 1. Existence of an approximate solution and a priori $H^1$ estimate.} Let us first observe that due to Lemma \ref{non_lin_curl},
 we have
	\begin{align}\label{form_adj_vw}
	-\left( \mathbf{curl}\,\sigma(\mathbf y)\times  \mathbf p_m,e_j\right)
	+\left(\mathbf{curl}\left(\sigma\left(\mathbf y\times \mathbf p_m\right)\right),e_j\right)&
	=b\left( p_m,e_j, \sigma(y) \right)-b\left(e_j,p_m, \sigma(y) \right)\nonumber\\
	&+
	b\left(p_m, y, \sigma(e_j) \right)-b\left( y,p_m, \sigma(e_j) \right)\end{align}
Let then $m$ be fixed and consider $Q: \ \mathbb{R}^m \longrightarrow \mathbb{R}^m$ defined by
	$$\begin{array}{ll}\left(Q \zeta\right)_i=&2\nu\left(Dp_m,De_j\right)+
	b\left( p_m,e_j, \sigma(y) \right)-b\left(e_j,p_m, \sigma(y) \right)\nonumber\\
	&+
	b\left(p_m, y, \sigma(e_j) \right)-b\left( y,p_m, \sigma(e_j) \right)-
	\left(f,  e_j\right),	\end{array}$$
where $ p_m=\sum_{i=1}^m\zeta_i   e_i$. The mapping $Q$ is obviously continuous. Let us prove that $Q(\zeta)\cdot \zeta>0$ if $|\zeta|$ is sufficiently large. Arguing as in the proof of  Proposition \ref{ex_uniq_lin}, we may prove that 
	\begin{align}\label{Q_zeta}{Q}(\zeta)\cdot \zeta&=
	2\nu \left\|D p_m\right\|_2^2+b\left(p_m, y, \sigma(p_m) \right)-b\left( y,p_m, \sigma(p_m) \right)-
	\left(f, p_m\right)\\
	&\geq \left(2\nu-\tfrac{2 \kappa_\star}{\nu} \left(\left\| u\right\|_2+\alpha
	\left\| \mathrm{curl} \, u\right\|_2\right)
	\right)\tfrac{\lambda_1}{2}\left|\zeta\right|^2-
	\left\|f\right\|_2\left|\zeta\right|\nonumber\\
	&\longrightarrow +\infty \qquad \mbox{when} \ 
	|\zeta|\rightarrow +\infty.\nonumber\end{align}
Due to the Brouwer theorem, we deduce that there exists $ \zeta^\ast\in \mathbb{R}^m$ such that ${Q}\left( \zeta^\ast\right)=0$ and thus $ p_m=\sum_{i=1}^m\zeta_i^\ast   e_i$ is a solution of problem (\ref{faedo_galerkin_adj}). Due to $(\ref{Q_zeta})$ and  Lemma \ref{rm2}, we deduce that
	$$\begin{array}{ll}2\nu \left\|D  p_m\right\|_2^2&=
	-\left( \mathbf{curl}\,\sigma(\mathbf p_m)\times  \mathbf y, \mathbf  p_m\right)+
	\left(f, p_m\right)\vspace{2mm}\\
	&\leq\tfrac{ 2\kappa_\star}{\nu}\left(\left\| u\right\|_2+\alpha
	\left\| \mathrm{curl} \, u\right\|_2\right)
	\left\|D  p_m\right\|_2^2+
	\kappa\left\|f\right\|_2
	 \left\|D p_m\right\|_2
	\end{array}$$
which gives
	\begin{equation}\label{est_H1_pm}\left(
	1-\tfrac{ \kappa_\star}{\nu^2}\left(\left\| u\right\|_2+\alpha
	\left\| \mathrm{curl} \, u\right\|_2\right)\right)\left\|Dp_m\right\|_2
	\leq \tfrac{\kappa}{2\nu} \left\|f\right\|_2.\end{equation}
{\it Step 2. A priori $H^2$ estimate.} Multiplying (\ref{faedo_galerkin_adj})  by $\left(1+\alpha\lambda_j\right) e_j$ yields
	$$-\nu\left(\Delta p_m,\mathbb{P}\sigma\left(e_j\right)\right)-
	\left( \mathbf{curl}\, \sigma(\mathbf y)\times  \mathbf p_m, \mathbb{P}\sigma\left(\mathbf e_j\right)\right)+
	\left( \mathbf{curl}\left(\sigma(\mathbf y\times \mathbf p_m\right), \mathbb{P}\sigma\left(\mathbf e_j\right)\right)=
	\left(f, \mathbb{P}\sigma\left(e_j\right)\right).$$
Multiplying by $\zeta^\ast_j$ and summing, it follows that
	$$-\nu\left(\Delta p_m,\mathbb{P}\sigma\left(p_m\right)\right)-
	\left( \mathbf{curl}\, \sigma(\mathbf y)\times  \mathbf p_m, \mathbb{P}\sigma\left(\mathbf p_m\right)\right)+
	\left( \mathbf{curl}\left(\sigma(\mathbf y\times \mathbf p_m\right), \mathbb{P}\sigma\left(\mathbf p_m\right)\right)=
	\left(f, \mathbb{P}\sigma\left(p_m\right)\right).$$
and thus 
	\begin{align}\label{adj_est_h2_1_1_adj}\left\| \sigma(p_m)\right\|_2^2&=
	\left(p_m, \sigma(p_m)\right)-\alpha\left(\Delta p_m, \sigma(p_m)\right)\nonumber\\
	&=\left(p_m, \sigma(p_m)\right)-\alpha\left(\Delta p_m, \sigma(p_m)-\mathbb{P} \sigma(p_m)\right)\nonumber\\
	&\quad +\tfrac{\alpha}{\nu}\left(\left(f,\mathbb{P} 
	\sigma(p_m)\right)-
	\left(\mathbf{curl}\,\sigma(\mathbf y)\times \mathbf p_m-
	\mathbf{curl}\left(\sigma\left(\mathbf y\times \mathbf p_m\right)\right),\mathbb{P} \sigma(\mathbf p_m)\right)\right). \end{align}
Arguing as in (\ref{est_h2_2_1_lin}) and (\ref{est_h2_4_1_lin} ), we obtain
\begin{equation}\label{est_h2_2_1_adj}
	\left|\left(p_m, \sigma(p_m)\right)-
	\alpha\left(\Delta p_m, \sigma(p_m)-\mathbb{P} \sigma(p_m)
	\right)\right|\leq c(1+\alpha)\left\|Dp_m\right\|_2\left\|\sigma(p_m)\right\|_2\end{equation}
and
	\begin{equation}\label{est_h2_3_1_adj} \left|\left(\mathbf{curl}\, \sigma(\mathbf y)\times  \mathbf p_m,\mathbb{P} 
	\sigma(\mathbf p_m)\right)\right|\leq\tfrac{c}{\alpha}
	\left\| \mathrm{curl}\, \sigma(y)\right\|_2
	\left\|\sigma(p_m)\right\|_2^2.\end{equation}
On the other hand, observing that
	$$b\left(y,\mathbb{P} \sigma(p_m),\sigma(p_m)\right)=b\left(y,\mathbb{P} \sigma(p_m)-\sigma(p_m),\sigma(p_m)\right)$$
and taking into account  Lemma \ref{prop_a_adj_2}, we obtain
	$$\begin{array}{ll}\left(\mathbf{curl}\left(\sigma\left(\mathbf y\times \mathbf p_m\right)\right),\mathbb{P} \sigma(\mathbf p_m)\right)&=
	b\left(\sigma(p_m),y,\mathbb{P} \sigma(p_m)\right)+b\left(y,\mathbb{P} \sigma(p_m)-\sigma(p_m),\sigma(p_m)\right)\vspace{2mm}\\
	&-b\left(\sigma(y),p_m,\mathbb{P} \sigma(p_m)\right)+b\left(p_m,\sigma(y),\mathbb{P} \sigma(p_m)\right)\vspace{2mm}\\
	&+ b\left(y, p_m,\mathbb{P} \sigma(p_m)\right)-b\left(p_m,y,\mathbb{P} \sigma(p_m)\right)\vspace{1mm}\\
	&-2\alpha\displaystyle\sum_{i=1}^2\left(b\left(
	\tfrac{\partial p_m}{\partial x_i},\tfrac{\partial y}{\partial x_i},\mathbb{P} \sigma(p_m)\right)
	-b\left(\tfrac{\partial y}{\partial x_i},\tfrac{\partial p_m}{\partial x_i},\mathbb{P} \sigma(p_m)\right)\right).\end{array}$$
Therefore, due to (\ref{sigma_psigma2})-(\ref{y2_sigma}), we have
	\begin{align}\label{ptimesy}&\left|\left(\mathbf{curl}\left(\sigma\left(\mathbf y\times \mathbf p_m\right)\right),
	\mathbb{P} \sigma(\mathbf p_m)\right)\right|\nonumber\\
	&\leq 	\left\|\sigma(p_m)\right\|_2\left\|\nabla y\right\|_\infty\left\|\mathbb{P} \sigma(p_m)\right\|_2+\left\|y\right\|_\infty
	\left\|\mathbb{P} \sigma(p_m)-\sigma(p_m)\right\|_{H^1}\left\|\sigma(p_m)\right\|_2\nonumber\\
	&+\left\|\sigma(y)\right\|_4\left\|\nabla p_m\right\|_4\left\|\mathbb{P} \sigma(p_m)\right\|_2+\left\|p_m\right\|_\infty
	\left\|\nabla\sigma(y)\right\|_2\left\|\mathbb{P} \sigma(p_m)\right\|_2\nonumber\\
	&+\left\|y\right\|_\infty\left\|\nabla p_m\right\|_2\left\|\mathbb{P} \sigma(p_m)\right\|_2+\left\|p_m\right\|_2
	\left\|\nabla y\right\|_\infty\left\|\mathbb{P} \sigma(p_m)\right\|_2\nonumber\\
	&+2\alpha\displaystyle\sum_{i=1}^2\left(\left\|
	\tfrac{\partial p_m}{\partial x_i}\right\|_4\left\|\nabla\left(\tfrac{\partial y}{\partial x_i}\right)\right\|_4
	\left\|\mathbb{P} \sigma(p_m)\right\|_2+
	\left\|\tfrac{\partial y}{\partial x_i}\right\|_\infty\left\|\nabla \left(\tfrac{\partial p_m}{\partial x_i}\right)\right\|_2
	\left\|	\mathbb{P} \sigma(p_m)\right\|_2\right)\nonumber\\
	&\leq c\left(\left\|y\right\|_{1,\infty}+\tfrac{1}{\alpha}
	\left\|\sigma(y)\right\|_{H^1}+\left\|y\right\|_{H^3}\right)\left\|\sigma(p_m)\right\|_2^2\nonumber\\
	&\leq \tfrac{c}{\alpha}\left(\left\|\mathrm{curl}\, \sigma(y)\right\|_{2}+\left\|Dy\right\|_2\right)
	\left\|\sigma(p_m)\right\|_2^2.
	\end{align}
Combining (\ref{adj_est_h2_1_1_adj})-(\ref{ptimesy}) and arguing as in the proof of (\ref{est_V2_zm}), we obtain
	\begin{equation}\label{est_V2_pm}
	\left(1-\tfrac{ \kappa_5}{\nu^2}\left(\left\| u\right\|_2+\alpha
	\left\| \mathrm{curl} \, u\right\|_2\right)\right)\left\| \sigma(p_m)\right\|_2\leq \kappa\left(\left(1+\alpha
	\right)
	\left\|Dp_m\right\|_2
	+\tfrac{\alpha}{\nu}\|f\|_2\right),\end{equation}
where $ \kappa_5$ and $\kappa$ are positive constants depending only on 
$\Omega$.\vspace{2mm}\\
{\it Step 3. Passing to the limit.} It remains to pass to the limit with respect to $m$. From estimates (\ref{est_H1_pm}) and (\ref{est_V2_pm}), it follows that if the condition 
	$$ \kappa_{\star\star}
	\left(\left\| u\right\|_2+\alpha
	\left\| \mathrm{curl} \, u\right\|_2\right)<\nu^2 \qquad \mbox{with} \ 
	 \kappa_{\star\star}=\max\left( \kappa_\star, \kappa_5\right)$$
is fulfilled, then there exists a subsequence, still indexed by $m$, and a function $ p\in  W$ such that
	$$ p_m \longrightarrow  p \qquad \mbox{weakly in} \  W.$$
By taking into account (\ref{form_adj_vw}) and passing to the limit in (\ref{faedo_galerkin_adj}), we obtain for every $j\geq 1$
	$$2\nu\left(Dp,D   e_j\right)+b\left( p,e_j, \sigma(y) \right)-b\left(e_j,p, \sigma(y) \right)+
	b\left(p, y, \sigma(e_j) \right)-b\left( y,p, \sigma(e_j) \right))=
	\left(f,  e_j\right)	$$
and by density we prove that $ p$ satisfies the variational formulation (\ref{form_var_lin_adj}). Moreover, $p$ satisfies estimates (\ref{est_H1_p}) and (\ref{est_V2_p}). Finally, since (\ref{adjoint}) is linear, the uniquess result is direct consequence of estimate (\ref{est_H1_p}). $\hfill\Box$

\section{Proof of the main results}
\setcounter{equation}{0}
{\bf Proof of the existence of an optimal control for $(P_\alpha)$.}\vspace{2mm}\\
We first prove Theorem \ref{main_existence}. Let $(u_{\alpha,k},y_{\alpha,k})_k\subset U_{ad}\times \left(W\cap H^3(\Omega)\right)$ be a minimizing sequence. Since $(u_{\alpha,k})_k$ is uniformly bounded in the closed convex set $U_{ad}$, we may extract a subsequence, still indexed by $k$, weakly convergent to some $u_\alpha\in U_{ad}$ in $H({\mathrm{curl}})$. On the other hand, due to estimate (\ref{state_est3}), we have
	$$\left\|y_{\alpha,k}\right\|_{H^3}\leq \tfrac{\kappa}{\alpha \nu}
	\left(\left\|u_{\alpha,k}\right\|_2+\alpha 
	\left\|\mathrm{curl}\, u_{\alpha,k}\right\|_2\right)$$
and the sequence $(y_{\alpha,k})_k$ is then bounded in $W\cap H^3(\Omega)$. Then there exists a subsequence, still indexed by $k$, and $y_\alpha\in W\cap H^3(\Omega)$ such that 
	$\left(y_{\alpha,k}\right)_{k}$ weakly converges to
$y_\alpha$ in $H^3(\Omega)$ and (by using compactness results on Sobolev spaces) strongly in $H^1(\Omega)$. 
Taking into account these convergence results and passing to the limit in the variational formulation corresponding to $y_{\alpha,k}$, we obtain
	$$2\nu \left(Dy_\alpha,D\phi\right)+
	b\left( \phi, y_\alpha, \sigma(y_\alpha)\right)-b\left(y_\alpha,
	 \phi, \sigma(y_\alpha)\right)=(u_\alpha,\phi) \qquad \mbox{for all} \ 
	\phi \in V$$
implying that $(u_\alpha,y_\alpha)$ satisfies (\ref{equation_etat}).
From the convexity and continuity of $J$, it follows the lower semicontinuity of $J$ in the weak topology and 
	$$J(u_\alpha,y_\alpha)\leq \liminf_kJ(u_{\alpha,k},y_{\alpha,k})=\inf(P_\alpha),$$
showing that $(u_\alpha,y_\alpha)$ is a solution for $(P_\alpha)$.$\hfill\Box$ \vspace{2mm}\\
{\bf Proof of the necessary optimality conditions  for $(P_\alpha)$.}\vspace{2mm}\\
Let us now prove Theorem \ref{main_1}. Assume that $\bar u_\alpha$ satisfies (\ref{control_constraint}) with $\bar\kappa=\max( \kappa_\star,\kappa^\star, \kappa_{\star\star},\kappa^{\star\star})$. Taking into account Proposition \ref{uniqueness_state}, Proposition \ref{ex_uniq_lin},   Proposition \ref{taylor} and Proposition \ref{ex_uniq_adj}, it follows that the corresponding state, linerized state, adjoint state exist and are unique and that the control-to-state mapping is G\^ateaux differentiable at $\bar u_\alpha$. For $\rho\in ]0,1[$ and $v\in U_{ad}$, let $u_{\alpha,\rho}=\bar u_\alpha+\rho (v-\bar u_\alpha)$, $y_{\alpha,\rho}$ the corresponding solution of (\ref{equation_etat}) and $z_{\alpha,\rho}=\tfrac{y_{\alpha,\rho}-\bar y_\alpha}{\rho}$. Since $(\bar u_\alpha,\bar y_\alpha)$ is an optimal solution and $(u_{\alpha,\rho},y_{\alpha,\rho})$ is admissible, we have 
	$$\displaystyle\lim_{\rho\rightarrow 0}
	\tfrac{J(u_{\alpha,\rho},y_{\alpha,\rho})-
	J(\bar u_\alpha,\bar y_\alpha)}{\rho}\geq 0.$$
By taking into account Proposition \ref{taylor}, we deduce that
	\begin{equation}\label{zv}\left(\bar z_{\alpha, v},\bar y_\alpha-y_d\right)+
	\lambda\left(\bar u_\alpha,v-\bar u_\alpha\right)\geq 0,\end{equation}
where $\bar z_{\alpha,v}$ is the (unique) solution of the linearized equation
	$$\left\{
  \begin{array}{ll}
    -\nu\Delta  \mathbf z+ \mathbf{curl}\,\sigma(\mathbf z)\times \bar{\mathbf y}_\alpha+
	 \mathbf{curl}\,\sigma(\bar{\mathbf y}_\alpha)\times \mathbf z+\nabla \pi=\mathbf v-\bar{\mathbf u}_\alpha&\quad\mbox{in} \ \Omega,\vspace{2mm}\\
	 \mathrm{div} \,  \mathbf z=0&\quad\mbox{in} \ \Omega,\vspace{2mm}\\
    \mathbf z\cdot \mathbf n=0, \qquad 
	\left(\mathbf n\cdot D\mathbf z\right)\cdot \boldsymbol\tau=0&\quad\mbox{on}\ \Gamma,
  \end{array}
\right .$$
Let then $\bar{p}_\alpha\in W$ be the unique solution of (\ref{adj_opt_eq_alpha}).
 Setting $\phi=\bar z_{\alpha, v}$ in the variational formulation 
(\ref{form_var_lin_adj}) and taking into account the variational formulation (\ref{var_lin}), we obtain
\begin{align}\left(\bar{y}_\alpha-y_d,\bar z_{\alpha,v}\right)
&=2\nu\left(D\bar p_\alpha,D\bar z_{\alpha,v}\right)+
	b\left(\bar p_\alpha,\bar z_{\alpha,v}, \sigma(\bar y_\alpha) \right)
	-b\left(\bar z_{\alpha,v},{\bar p}_\alpha, \sigma(\bar y_\alpha) \right)\nonumber\\
& \ +
	b\left(\bar p_\alpha, \bar y_\alpha, \sigma(\bar z_{\alpha,v}) \right)
	-b\left(\bar y_\alpha,\bar p_\alpha, \sigma(\bar z_{\alpha,v}) \right)\nonumber\\
&\label{GF_op}=\left(v-\bar{u}_\alpha,\bar{p}_\alpha\right).
	\end{align}
The result follows by combining (\ref{zv}) and (\ref{GF_op}).\vspace{2mm}\\
Let us finally prove the results related with the asymptotic analysis when $\alpha$ tends to zero.\vspace{2mm}\\
{\bf Proof of Theorem \ref{assympt_1}.}\vspace{2mm}\\
 The proof is split into two steps.\vspace{2mm}\\
{\it Step 1. Convergence of $(\bar u_\alpha,\bar y_\alpha)$.}
By taking into account (\ref{state_est1}) and (\ref{state_est2}), we have
	$$\left\|D\bar y_\alpha\right\|_2\leq \tfrac{ \kappa_2}{\nu} \left\|\bar u_\alpha\right\|_2, $$
	$$ \left\|\mathrm{curl}\, \sigma(\bar y_\alpha)
	\right\|_2\leq \tfrac{\kappa}{\nu}\left(
	\left\|\bar u_\alpha\right\|_2+\alpha \left\|\mathrm{curl}\,\bar u_\alpha\right\|_2\right),$$
and thus $\left(\bar y_\alpha\right)_\alpha$ and
$\left(\mathrm{curl}\, \sigma(\bar y_\alpha)\right)_\alpha$ are bounded independently of $\alpha$ (when $\alpha$ tends to zero). There then exists a subsequence, still indexed by $\alpha$, $\bar u_0\in U_{ad}$, $\bar y_0\in V$ and $\bar\omega_0\in L^2(\Omega)$ such that
	$$\bar u_\alpha\longrightarrow 
	\bar u_0 \qquad
	 \mbox{weakly in} \ L^2(\Omega). $$
	$$\bar y_\alpha \longrightarrow \bar y_0 \qquad
	 \mbox{weakly in} \ H^1(\Omega) \
	\mbox{and strongly in} \ L^2(\Omega),$$
	$$\mathrm{curl}\, \sigma(\bar y_\alpha)\longrightarrow 
	\bar\omega_0 \qquad
	 \mbox{weakly in} \ L^2(\Omega). $$
By taking into account $(\ref{var_form_state})$ and $(\ref{transport_state})$, we have
	\begin{equation}\label{vf_yalpha}2\nu\left(D\bar y_\alpha,D\phi\right)+\left(\mathbf{curl}\, \sigma(\mathbf{\bar y_\alpha})\times \mathbf{\bar y_\alpha},\phi\right)=
	\left(\bar u_\alpha,\phi\right) \qquad \mbox{for all} \ 
	\phi\in V\end{equation}
and 
	\begin{equation}\label{vf_sigma_alpha}\left(\mathrm{curl}\, \sigma(\bar y_\alpha),\phi\right)-
	\tfrac{\alpha}{\nu}\, b(\bar y_\alpha,\phi,\mathrm{curl}\, \sigma(\bar y_\alpha))=\left(\tfrac{\alpha}{\nu}\,\mathrm{curl}\, \bar u_\alpha+
	\mathrm{curl}\, \bar y_\alpha,\phi\right)\qquad 
	\mbox{for all} \ 
	\phi\in {\cal D}(\Omega).\end{equation}
The previous convergence results yield
	$$\lim_{\alpha\rightarrow 0^+}\left(\mathbf{curl}\, \sigma(\mathbf{\bar y_\alpha})\times \mathbf{\bar y_\alpha},\phi\right)=\left(\boldsymbol{\bar\omega}_0\times \mathbf{\bar y}_0,\phi\right)\qquad \mbox{for all} \ 
	\phi\in V$$
and 
	$$\lim_{\alpha\rightarrow 0^+} 
	b(\bar y_\alpha,\phi,\mathrm{curl}\, \sigma(\bar y_\alpha))=
	b(\bar y_0,\phi,\bar \omega_0)\qquad 
	\mbox{for all} \ 
	\phi\in {\cal D}(\Omega).$$
By passing to the limit in (\ref{vf_yalpha}) and (\ref{vf_sigma_alpha}), we deduce that
	$$2\nu\left(D\bar y_0,D\phi\right)+\left(\boldsymbol{\bar\omega}_0\times
	 \mathbf{\bar y}_0,\phi\right)=
	\left(\bar u_0,\phi\right) \qquad \mbox{for all} \ 
	\phi\in V$$
and 
	$$\left(\bar\omega_0,\phi\right)=\left(\mathrm{curl}\, 
	\bar y_0,\phi\right)\qquad 
	\mbox{for all} \ 
	\phi\in {\cal D}(\Omega).$$
Therefore, $\bar\omega_0=\mathrm{curl}\, 
	\bar y_0$ and $\bar y_0$ satisfies
	$$2\nu\left(D\bar y_0,D\phi\right)+b(\bar y_0,\bar y_0,\phi)=
	\left(\bar u_0,\phi\right) \qquad \mbox{for all} \ 
	\phi\in V$$ 
that is, $(\bar u_0,\bar y_0)$ is admissible for $(P_0)$. Let us now prove that the convergence of $\bar y_\alpha$ to $\bar y_0$ is strong. Taking into account the variational formulations corresponding to $\bar y_\alpha$ and $\bar y_0$, we easily see that $\bar y_\alpha-\bar y_0$ satisfies
	$$\begin{array}{ll}2\nu \left\|D\left(\bar y_\alpha-\bar y_0\right)\right\|_2^2
	&=\left(\bar u_\alpha-\bar u_0,\bar y_\alpha-\bar y_0\right)-
	\left(\mathrm{curl}\, \sigma(\bar y_\alpha)\times \bar y_0
	-\mathrm{curl}\, \bar y_0\times \bar y_0,
	\bar y_\alpha-\bar y_0\right)\vspace{2mm}\\
	&\longrightarrow 0 \qquad \mbox{when} \ \alpha \rightarrow 0^+.
	\end{array}$$
{\it Step 2. Convergence of $\bar p_\alpha$.} 
By taking into account (\ref{est_H1_p}), we have
	$$\left(1-\tfrac{ \kappa_{\star\star}}{\nu^2}\left(
	\left\|\bar u_\alpha\right\|_2+\alpha\left\|{\rm curl}\, \bar u_\alpha\right\|_2\right)\right) \left\|D\bar p_\alpha\right\|_2\leq
	 \tfrac{\kappa}{\nu}\left\|\bar y_\alpha-y_d\right\|_2,$$
and thus $\left(\bar p_\alpha\right)_\alpha$ is also bounded in $H^1(\Omega)$ independently of $\alpha$ (when $\alpha$ tends to zero). There then exists a subsequence, still indexed by $\alpha$, and $\bar p_0\in V$ such that
	$$\bar p_\alpha \longrightarrow \bar{p}_0 \qquad
	 \mbox{weakly in} \ H^1(\Omega) \
	\mbox{and strongly in} \ L^2(\Omega).$$
By taking into account the convergence results established in the first step, we deduce that
	$$\begin{array}{ll}\displaystyle\lim_{\alpha\rightarrow 0^+}\left(b\left( \bar p_\alpha,\phi, \sigma(\bar y_\alpha) \right)-b\left( \phi,\bar p_\alpha, \sigma(\bar y_\alpha) \right)\right)&=\displaystyle\lim_{\alpha\rightarrow 0^+}\left(\mathrm{curl}\,\sigma(\bar y_\alpha) \times \bar p_\alpha,\phi\right)\vspace{2mm}\\
	&=\displaystyle\lim_{\alpha\rightarrow 0^+} \left(\mathrm{curl}\,\bar y_0\times \bar p_0,\phi\right)\vspace{2mm}\\
	&=b(\bar p_0,\bar y_0,\phi)-b(\phi,\bar y_0,\bar p_0)
		\end{array}$$
and 
	$$\lim_{\alpha\rightarrow 0^+}\left(b\left(\bar p_\alpha, \bar y_\alpha, \sigma(\phi) \right)-b\left( \bar y_\alpha,\bar p_\alpha, \sigma(\phi) \right)\right)=b\left(\bar p_0, \bar y_0,\phi \right)-b\left(\bar  y_0,\bar p_0,\phi\right)$$
for all $\phi\in W$. Passing then to the limit in the variational formulation (\ref{form_var_lin_adj}) corresponding to $\bar p_\alpha$ yields
$$
	2\nu\left(D \bar p_0,D  \phi\right)+b(\phi,\bar y_0,\bar p_0)-b\left(\bar y_0,\bar p_0,\phi\right)
	=\left(\bar y_0-y_d, \phi\right)$$
for all $\phi\in W$ and thus $\bar p_0$ is the unique weak solution of 
	(\ref{adjoint_limit}). The optimality condition for the control follows by passing to the limit in $(\ref{opt_control_alpha})$.$\hfill\Box$
\begin{remark} Notice that if $\bar u_\alpha$ satisfies $(\ref{control_constraint})$, then the limit $\bar u_0$ satisfies 
	$\bar\kappa\, \left\|\bar u_0\right\|_2<\nu^2$.
Taking into account Remark $\ref{remark1}$, we deduce that
	$\tfrac{ \kappa_1 \kappa_2}{2}\, \left\|\bar u_0\right\|_2<\nu^2$,
which implies the uniqueness of $\bar y_0$ and $\bar p_0$.\end{remark}
{\bf Proof of Corollary \ref{corollary_assymp}.}\vspace{2mm}\\ By taking into account the convergence results of Theorem  $\ref{assympt_1}$ and the lower semicontinuity of $J$, we obtain
	$$\min(P_0)\leq J(\bar u_0,\bar y_0)\leq \liminf_{\alpha\rightarrow 0^+}J(\bar u_\alpha,\bar y_\alpha)=\liminf_{\alpha\rightarrow 0^+}\min(P_\alpha).$$ 
On the other hand, let $(\hat u,\hat y)$ be a solution of problem $(P_0)$ and let $\hat y_\alpha$ be a solution of $(\ref{equation_etat})$ corresponding to $\hat u$. Then $(\hat u,\hat y_\alpha)$ is admissible for $(P_\alpha)$ and 
	\begin{equation}\label{semi_con_1}
	\min(P_\alpha)\leq J(\hat u,\hat y_\alpha).\end{equation}
Arguing as in the first step of the proof of Theorem \ref{assympt_1}, we can establish the convergence of $\hat y_\alpha$ to
 $\hat y$ in $V$ and thus 
	\begin{equation}\label{semi_con_2}\lim_{\alpha\rightarrow 0^+}\min(P_\alpha)\leq \lim_{\alpha\rightarrow 0^+}J(\hat u,\hat y_\alpha)=J(\hat u,\hat y)=\min(P_0).\end{equation}
The conclusion follows from (\ref{semi_con_1}) and (\ref{semi_con_2}).$\hfill\Box$

\appendix
\section{Appendix} \label{appendix}
\setcounter{equation}{0}
The aim of this section is to prove Lemma \ref{non_lin_curl}, Lemma \ref{prop_a_adj_2} and Lemma \ref{rm2}. \vspace{2mm}\\
{\bf Proof of Lemma \ref{non_lin_curl}.} Let us first observe that 
	$$\begin{array}{ll}\mathbf{curl}\left(\mathbf z\times \boldsymbol\phi\right)&=\left( \mathrm{div}\, \boldsymbol\phi\right)\mathbf z
	+\boldsymbol\phi\cdot \nabla \mathbf z
	-\left( \mathrm{div}\, \mathbf z\right)\boldsymbol\phi-\mathbf z\cdot \nabla \boldsymbol\phi\vspace{2mm}\\
	&=
	\boldsymbol\phi\cdot \nabla \mathbf z-\mathbf z\cdot \nabla \boldsymbol\phi\end{array}$$
and 
	$$\left(\mathbf{curl}\,\mathbf w,\mathbf u\right)=\left(\mathbf w,\mathbf{curl}\,\mathbf u\right)+\int_\Gamma
	\left(\mathbf w\times \mathbf u\right)\cdot \mathbf n\,dS.$$
It follows that
	$$\begin{array}{ll}\left(\mathbf{curl}\,\sigma(\mathbf y)\times \mathbf z,\boldsymbol\phi\right)&=\left(\mathbf{curl}\,\sigma(\mathbf y),\mathbf  z \times \boldsymbol\phi\right)\vspace{2mm}\\
	&=\displaystyle \left(\sigma(\mathbf y),\mathbf{curl}\left(\mathbf  z \times \boldsymbol\phi\right)\right)+\int_\Gamma
	\sigma(\mathbf y)\times \left(\mathbf  z \times \boldsymbol\phi\right)\cdot \mathbf n\,dS\vspace{1mm}\\
	&= \left(\sigma(\mathbf y),\mathbf{curl}\left(\mathbf  z \times \boldsymbol\phi\right)\right)\vspace{2mm}\\
	&= \left(\sigma(\mathbf y),\boldsymbol\phi\cdot \nabla \mathbf z-\mathbf z\cdot \nabla \boldsymbol\phi\right)\vspace{2mm}\\
	&=b(\boldsymbol\phi,\mathbf z,\sigma(\mathbf y))-
	b(\mathbf z,\boldsymbol\phi,\sigma(\mathbf y))\end{array}$$
which gives the first identity. The proof of the second identity is split into two steps. 
In order to give a sense to the different boundary terms, we 
 first assume that $y, z$ belong to $W\cap H^4(\Omega)$. We next apply a regularization process to prove that the results are still valid for $y, z\in W\cap H^3(\Omega)$.\vspace{1mm}\\
{\it Step 1.} Arguing as above, we obtain
	\begin{align}\label{curl_sigma}\left( \mathbf{curl}\, \sigma\left(\mathbf y\times \mathbf z\right), \boldsymbol\phi\right)&=
	\left( \mathbf{curl}\, \left(\mathbf y\times \mathbf z\right), \boldsymbol\phi\right)-\alpha
	\left( \mathbf{curl}\, \Delta\left(\mathbf y\times \mathbf z\right), \boldsymbol\phi\right)\nonumber\\
	&=b\left(z,y,\phi\right)-b\left(y,z,\phi\right)-\alpha
	\left( \mathbf{curl}\, \Delta\left(\mathbf y\times \mathbf z\right), \boldsymbol\phi\right),\end{align}
where
	\begin{align}\label{curl_delta}\left(\mathbf{curl}\left(\Delta\left(\mathbf y\times \mathbf z\right)\right), \boldsymbol\phi\right)&=\displaystyle\left(\Delta\left(\mathbf y\times \mathbf z\right),\mathbf{curl}\, \boldsymbol\phi\right)+I_1\nonumber\\
	&=\left(-\mathbf{curl}\left(\mathbf{curl}\left(\mathbf y\times \mathbf z\right)\right)
	+\nabla\left( \mathrm{div}\left(\mathbf y\times \mathbf z\right)\right),\mathbf{curl}\, \boldsymbol\phi\right)+I_1\nonumber\\
	&=\left(-\mathbf{curl}\left(\mathbf{curl}\left(\mathbf y\times \mathbf z\right)\right),\mathbf{curl}\,
	 \boldsymbol\phi\right)+I_1\nonumber\\
	&=-\left(\mathbf{curl}\left(\mathbf y\times \mathbf z\right),\mathbf{curl}\left(\mathbf{curl}\, \boldsymbol\phi\right)\right)
	\displaystyle -I_2+I_1\nonumber\\
	&=\left(\mathbf{curl}\left(\mathbf y\times \mathbf z\right),\Delta \boldsymbol\phi-\nabla\left( \mathrm{div}\, \boldsymbol\phi\right)\right)+I_1-I_2\nonumber\\
	&=b(\mathbf z,\mathbf y,\Delta\boldsymbol\phi)-
	b(\mathbf y,\mathbf z,\Delta \boldsymbol\phi)+I_1-I_2\nonumber\\
	&=b(z,y,\Delta\phi)-b(y,z,\Delta\phi)+I_1-I_2
	\end{align}
with  $$\begin{array}{ll}I_1=\displaystyle \int_\Gamma
	\left(\Delta\left(\mathbf y\times \mathbf z\right)\times \boldsymbol\phi\right)\cdot \mathbf n\,dS
	\qquad\mbox{and}\qquad 
	I_2=\displaystyle \int_\Gamma \mathbf{curl}\left(\mathbf y\times \mathbf z\right)\times \mathbf{curl}\,\boldsymbol\phi\cdot \mathbf n\,dS. \end{array}$$
Let us now prove that $I_1-I_2=0$.
 By taking into account Lemma \ref{curl_trace}, we have
	\begin{align}\label{I_2}
	\mathbf{curl}\left(\mathbf y\times \mathbf z\right)\times 
	\mathbf{curl}\, \boldsymbol \phi \cdot \mathbf n\big|_{\Gamma}
	&=\left(\mathbf z\cdot \nabla \mathbf y-\mathbf y\cdot \nabla \mathbf z\right)\times 
	(0,0,\phi\cdot g)^\top \cdot \mathbf n\big|_{\Gamma}\nonumber\\
	&=\left(\phi\cdot g\right)\left(\mathbf z\cdot \nabla \mathbf y-\mathbf y\cdot \nabla \mathbf z\right)\cdot \boldsymbol\tau\big|_{\Gamma}\nonumber\\
	&=\left(\phi\cdot g\right)
	\left(z\cdot \nabla y-y\cdot \nabla z\right)\cdot 
	\tau\big|_{\Gamma}.
	\end{align}
Similarly, since
$$\begin{array}{ll}\Delta\left(\mathbf y\times \mathbf z\right)&
	=-\mathbf{curl}\left(\mathbf{curl}\left(\mathbf y\times \mathbf z\right)\right)
	+\nabla\left( \mathrm{div}\left(\mathbf y\times \mathbf z\right)\right)\vspace{3mm}\\
	&=-\mathbf{curl}\left(\mathbf z\cdot \nabla \mathbf y-\mathbf y\cdot \nabla \mathbf z\right)\vspace{1mm}\\
	&=-\mathbf z\cdot \nabla \left(\mathbf{curl}\,\mathbf y\right)+\mathbf{curl}\,\mathbf y\cdot \nabla \mathbf z-
	\left( \mathrm{div}\,\mathbf z\right)\cdot \mathbf{curl} \,\mathbf y-\displaystyle \sum_{k=1}^2
	\nabla \mathbf z_k\times \nabla\mathbf y_k\vspace{1mm}\\
	&+\mathbf y\cdot \nabla \left(\mathbf{curl}\,\mathbf z\right)-\mathbf{curl}\,\mathbf z\cdot \nabla \mathbf y+
	\left( \mathrm{div}\,\mathbf y\right)\cdot \mathbf{curl} \,\mathbf z+\displaystyle \sum_{k=1}^2
	\nabla \mathbf y_k\times \nabla\mathbf z_k\vspace{1mm}\\
	&=\mathbf y\cdot \nabla \left(\mathbf{curl}\,\mathbf z\right)-\mathbf z\cdot \nabla \left(\mathbf{curl}\,\mathbf y\right)+2\displaystyle \sum_{k=1}^2
	\nabla \mathbf y_k\times \nabla\mathbf z_k\end{array}$$
we deduce that
	$$\begin{array}{ll}\Delta\left(\mathbf y\times \mathbf z\right)\big|_{\Gamma}-\displaystyle 2\sum_{k=1}^2\left(\nabla \mathbf y_k\times \nabla\mathbf z_k\right)\big|_{\Gamma}
	&=\mathbf y\cdot \nabla \left(\mathbf{curl}\,\mathbf z\right)\big|_{\Gamma}-\mathbf z\cdot \nabla \left(\mathbf{curl}\,\mathbf y\right)\big|_{\Gamma}\vspace{0mm}\\
	&=\left(0,0,y\cdot \nabla\left(z\cdot g\right)-z\cdot \nabla\left(y\cdot g\right)\right)^\top\big|_{\Gamma}\vspace{2mm}\\
	&=\left(0,0,\left(y\cdot \nabla z\right)\cdot g+\left(y\cdot \nabla g\right)\cdot z\right)^\top\big|_{\Gamma}\vspace{2mm}\\
	&-\left(0,0,-\left(z\cdot \nabla y\right)\cdot g-\left(y\cdot \nabla g\right)\cdot y\right)^\top\big|_{\Gamma}
	\end{array}
	$$
and thus
	\begin{align}\label{J1J2}&\left(\Delta\left(\mathbf y\times \mathbf z\right)\times \boldsymbol\phi\right)\cdot \mathbf n\big|_{\Gamma}-2\displaystyle \sum_{k=1}^2\left(\left(\nabla \mathbf y_k\times \nabla\mathbf z_k\right)\times \boldsymbol \phi\right)\cdot \mathbf n\big|_{\Gamma}\nonumber\\
	&=
	-\left(y\cdot \nabla z-z\cdot \nabla y\right)\cdot g \left(\phi\cdot \tau\right)\big|_{\Gamma}-\left(\left(y\cdot \nabla g\right)
	\cdot z-\left(z\cdot \nabla g\right)\cdot y\right)\left(\phi\cdot \tau\right)\big|_{\Gamma}.\end{align}
Easy calculations, together with the fact that $y\big|_{\Gamma}=
\left(y\cdot \tau\right)\tau\big|_{\Gamma}$ and $z\big|_{\Gamma}=
\left(z\cdot \tau\right)\tau\big|_{\Gamma}$,  show that
	\begin{align}\label{J1}
	\left(y\cdot \nabla g\right)\cdot z\big|_{\Gamma}-
	\left(z\cdot \nabla g\right)\cdot y\big|_{\Gamma}
	&=\left(\nabla g\, y\right)\cdot z\big|_{\Gamma}-
	\left(\nabla g\, z\right)\cdot y\big|_{\Gamma}\nonumber\\
	&=\left(\nabla g\, \left(y\cdot \tau\right)\tau\right)\cdot \left(z\cdot \tau\right)\tau\big|_{\Gamma}-
	\left(\nabla g\, \left(z\cdot \tau\right)\tau\right)\cdot \left(y\cdot \tau\right)\tau\big|_{\Gamma}
	\nonumber\\
	&=\left(y\cdot \tau\right)\left(z\cdot \tau\right)
	\left(\left(\nabla g\, \tau\right)\cdot \tau-
	\left(\nabla g\, \tau\right)\cdot \tau\right)\big|_{\Gamma}=0.
	\end{align}
Taking into account (\ref{I_2})-(\ref{J1}), we deduce that
	$$\begin{array}{ll}
	I_1-I_2&=\displaystyle 2\sum_{k=1}^2\int_\Gamma \left(\nabla \mathbf y_k\times \nabla \mathbf z_k\right)\times \boldsymbol\phi\cdot \mathbf n \,dS\vspace{2mm}\\
	&+\displaystyle\int_\Gamma
	 \left(\left(\phi\cdot \tau\right)\left(z\cdot \nabla y-y\cdot \nabla z\right)\cdot g-\left(\phi\cdot g\right)\left(z\cdot \nabla y-y\cdot \nabla z\right)\cdot \tau\right)dS.
	\end{array}$$
On the other hand, due to Lemma 4.1 in \cite{K06} we have
	$$\left(z\cdot \nabla y-y\cdot \nabla z\right)\cdot n\big |_\Gamma=0.$$
Observing that $u\cdot n\big|_{\Gamma}=0$ implies that
	$$\begin{array}{ll}\left(\phi\cdot \tau\right)\left(u\cdot g\right)
	 -\left(\phi\cdot g\right)\left(u\cdot \tau\right)\big|_{\Gamma}&
	=\left(\phi\cdot \tau\right)\left(\left(u\cdot \tau\right)\tau\cdot g\right)-\left(\left(\phi\cdot \tau\right)\tau\cdot g\right)\left(u\cdot \tau\right)\big|_{\Gamma}
	\vspace{2mm}\\
	&=\left(\phi\cdot \tau\right)\left(u\cdot \tau\right)
	\left(\tau\cdot g-\tau\cdot g\right)\big|_{\Gamma}=0,\end{array}$$
we deduce that
	$$I_1-I_2	=\displaystyle 2\sum_{k=1}^2\int_\Gamma 
	\left(\nabla \mathbf y_k\times \nabla \mathbf z_k\right)
	\times \boldsymbol\phi\cdot \mathbf n \,dS.$$
Finally, since $\left(Du\cdot n\right)\cdot \tau\big|_{\Gamma}=0$ and $\mathrm{div}\, u=0$ imply 
	$$\left(n_2^2-n_1^2\right)\left(\tfrac{\partial u_1}{\partial x_2}+
	\tfrac{\partial u_2}{\partial x_1}\right)-2n_1n_2\,\tfrac{\partial u_1}{\partial x_1}=0 \qquad \mbox{on} \ \Gamma,$$
we obtain
	$$\begin{array}{ll}
	\displaystyle \sum_{k=1}^2\left(\nabla \mathbf y_k\times \nabla\mathbf z_k\right)\times \boldsymbol \phi\cdot \mathbf n\big|_{\Gamma}
	&\displaystyle =\sum_{k=1}^2\left(\tfrac{\partial y_k}{\partial x_1}\tfrac{\partial z_k}{\partial x_2}-
	\tfrac{\partial y_k}{\partial x_2}\tfrac{\partial z_k}{\partial x_1}\right) \phi\cdot \tau\big|_{\Gamma}\vspace{2mm}\\
	&\displaystyle =\left(\tfrac{\partial y_1}{\partial x_1}\left(\tfrac{\partial z_1}{\partial x_2}+
	\tfrac{\partial z_2}{\partial x_1}\right)-\tfrac{\partial z_1}{\partial x_1}\left(\tfrac{\partial y_1}{\partial x_2}+
	\tfrac{\partial y_2}{\partial x_1}\right)\right) \phi\cdot \tau\big|_{\Gamma}\vspace{2mm}\\
	&=0
\end{array}$$
and thus
	\begin{equation}\label{I1-I2}I_1-I_2=0.\end{equation}
The conclusion is then a consequence of (\ref{curl_sigma}), (\ref{curl_delta}) and (\ref{I1-I2}).\vspace{1mm}\\
{\it Step 2. Regularization process.} Let us now go back to the case 
$y, z\in W\cap H^3(\Omega)$. We first infer that there exist $y_\varepsilon, z_\varepsilon \in W\cap H^4(\Omega)$ such that
	$$\lim_{\varepsilon\rightarrow 0^+}
	\left\|y-y_\varepsilon\right\|_{H^3}
	=\lim_{\varepsilon\rightarrow 0^+}
	\left\|z-z_\varepsilon\right\|_{H^3}=0.$$
Indeed, if $y\in W\cap H^3(\Omega)$, then it  satisfies the following Stokes system
	\begin{equation}\label{stokes_molifier}\left\{\begin{array}{ll}
   	 -\Delta y+y+\nabla \pi=f&\quad\mbox{in} \ \Omega,\vspace{2mm}\\
	 \mathrm{div} \,  y=0&\quad\mbox{in} \ \Omega,\vspace{2mm}\\
   	 y\cdot n=0, \qquad 
	\left(n\cdot D y\right)\cdot \tau=0&\quad\mbox{on}\ \Gamma
 	 \end{array}
	\right. \end{equation}
with $f=-\Delta y+y\in H^1(\Omega)$ and $\pi=0$. Using Friedrichs mollifiers, we can construct 
$f_\varepsilon\in H^2(\Omega)$ such that
	$$\lim_{\varepsilon \rightarrow 0^+}
	\left\|f-f_\varepsilon\right\|_{H^1}=0.$$
Let $y_\varepsilon\in H^4(\Omega)$ be the solution of (\ref{stokes_molifier}) corresponding to $f_\varepsilon$. By using classical regularity results, we have
	$$\left\|y-y_\varepsilon\right\|_{H^3} \leq c 
	\left\|f-f_\varepsilon\right\|_{H^1} \longrightarrow 0 \quad \mbox{when} \ \varepsilon \rightarrow 0.$$
On the other hand, by taking into account the first step, we deduce that
	$$\left( \mathbf{curl}\, \sigma\left(\mathbf y_\varepsilon\times \mathbf z_\varepsilon\right), \boldsymbol\phi\right)
	=b\left(z_\varepsilon,y_\varepsilon, \sigma(\phi)\right)-b\left(y_\varepsilon,z_\varepsilon,\sigma(\phi)\right)$$
and the result follows by passing to the limit.$\hfill\Box$ \vspace{2mm}\\
{\bf Proof of Lemma \ref{prop_a_adj_2}.} Standard calculations show that
	$$\Delta\left(\mathbf y\times \mathbf z\right)=\mathbf y\times \Delta \mathbf z-\mathbf z\times \Delta\mathbf y-
	2\sum_{i=1}^2\tfrac{\partial \mathbf z}{\partial x_i}\times \tfrac{\partial \mathbf y}{\partial x_i}$$
and thus 
	$$\begin{array}{ll}&\mathbf{curl}\left(\sigma\left(\mathbf y\times \mathbf z\right)\right)\vspace{0mm}\\
	&=\displaystyle
	\mathbf{curl}\left(\mathbf y\times \sigma(\mathbf z)\right)-
	\mathbf{curl}\left(\mathbf z\times \sigma(\mathbf y)\right)+\mathbf{curl}\left(\mathbf z\times \mathbf y\right)-2\alpha \displaystyle
	\sum_{i=1}^2\mathbf{curl}\left(\tfrac{\partial \mathbf z}{\partial x_i}\times \tfrac{\partial \mathbf y}{\partial x_i}\right)
	\vspace{1mm}\\
	&=\displaystyle \left( \mathrm{div}\,\sigma(\mathbf z)\right) \mathbf y+\sigma(\mathbf z)\cdot \nabla \mathbf y-\left( \mathrm{div}\,\mathbf y\right) \sigma(\mathbf z)-\mathbf y\cdot \nabla \left(\sigma(\mathbf z)\right)\vspace{2mm}\\
	&-\displaystyle \left( \mathrm{div}\,\sigma(\mathbf y)\right) \mathbf z-\sigma(\mathbf y)\cdot \nabla \mathbf z+\left( \mathrm{div}\,\mathbf z\right) \sigma(\mathbf y)+\mathbf z\cdot \nabla \left(\sigma(\mathbf y)\right)\vspace{2mm}\\
	&\displaystyle+ \left( \mathrm{div}\,\mathbf y\right) \mathbf z+\mathbf y\cdot \nabla \mathbf z-
\left( \mathrm{div}\,\mathbf z\right) \mathbf y-\mathbf z\cdot \nabla \mathbf y\vspace{2mm}\\
	&\displaystyle-2\alpha\displaystyle\sum_{i=1}^2\left( \mathrm{div}\left(\tfrac{\partial \mathbf z}{\partial x_i}\right)
	\tfrac{\partial \mathbf y}{\partial x_i}+\tfrac{\partial \mathbf z}{\partial x_i}\cdot \nabla \left(\tfrac{\partial \mathbf y}{\partial x_i}\right)
	- \mathrm{div}\left(\tfrac{\partial \mathbf y}{\partial x_i}\right)
	\tfrac{\partial \mathbf z}{\partial x_i}-\tfrac{\partial \mathbf y}{\partial x_i}\cdot \nabla \left(\tfrac{\partial \mathbf z}{\partial x_i}\right)\right)
	\vspace{2mm}\\
	&=\displaystyle \sigma(\mathbf z)\cdot \nabla \mathbf y-\mathbf y\cdot \nabla \left(\sigma(\mathbf z)\right)-\sigma(\mathbf y)\cdot \nabla \mathbf z+\mathbf z\cdot \nabla \left(\sigma(\mathbf y)\right)+ \mathbf y\cdot \nabla \mathbf z-\mathbf z\cdot \nabla \mathbf y\vspace{2mm}\\
	&\displaystyle-2\alpha\displaystyle\sum_{i=1}^2\left(
	\tfrac{\partial \mathbf z}{\partial x_i}\cdot \nabla \left(\tfrac{\partial \mathbf y}{\partial x_i}\right)
	-\tfrac{\partial \mathbf y}{\partial x_i}\cdot \nabla \left(\tfrac{\partial \mathbf z}{\partial x_i}\right)\right).\end{array}$$
Therefore
	$$\begin{array}{ll}\left(\mathbf{curl}\left(\sigma\left(\mathbf y\times \mathbf z\right)\right),\boldsymbol\phi\right)&=
	b\left(\sigma(z),y,\phi\right)-b\left(y,\sigma(z),\phi\right)-b\left(\sigma(y),z,\phi\right)+b\left(z,\sigma(y),\phi\right)\vspace{2mm}\\
	&+ b\left(y, z,\phi\right)-b\left(z,y,\phi\right)-2\alpha\displaystyle\sum_{i=1}^2\left(b\left(
	\tfrac{\partial z}{\partial x_i},\tfrac{\partial y}{\partial x_i},\phi\right)
	-b\left(\tfrac{\partial y}{\partial x_i},\tfrac{\partial z}{\partial x_i},\phi\right)\right)
	\end{array}$$
which gives the claimed result.$\hfill\Box$\vspace{2mm}\\
{\bf Proof of Lemma \ref{rm2}.}  
By taking into account Lemma \ref{non_lin_curl}, we obtain
	\begin{align}	\label{trilin_0}
	&\left( \mathbf{curl}\,\sigma(\mathbf z)\times  \mathbf y,\mathbf z\right)=b\left(z, y, \sigma(z)\right) 
	-b\left(y,z, \sigma(z)\right) 
	\notag\\
	&=b(z,y,z)-b(y,z,z)-\alpha 
	\left(z\cdot \nabla y-y\cdot \nabla z,\Delta z\right)
\notag\\
	&=b(z,y,z)+\alpha 
	\left(\mathbf z\cdot \nabla \mathbf y-\mathbf y\cdot \nabla \mathbf z,\mathbf{curl}\left(\mathbf{curl}\, \mathbf z\right)\right)
	\notag\\
&=b(z,y,z)+\alpha 
	\left(\mathbf{curl}\left(\mathbf z\cdot \nabla \mathbf y-\mathbf y\cdot \nabla \mathbf z\right),\mathbf{curl}\, \mathbf z\right)
	+\alpha I\notag\\
	&=b(z,y,z)
	+\alpha\left(b\left(\mathbf z, \mathbf{curl}\, \mathbf y,\mathbf{curl}\, \mathbf z\right)
	 -b\left(\mathbf y, \mathbf{curl}\, \mathbf z,\mathbf{curl}\, \mathbf z\right)\right)\notag\\
	&+2\alpha\sum_{k=1}^3\left(\nabla \mathbf z_k\times \nabla \mathbf y_k,\mathbf{curl}\,\mathbf z\right)+\alpha I\notag\\
	&=b(z,y,z)
	+\alpha b\left(\mathbf z, \mathbf{curl}\, \mathbf y,\mathbf{curl}\, \mathbf z\right)+2\alpha\sum_{k=1}^2\left(\nabla \mathbf z_k\times \nabla \mathbf y_k,\mathbf{curl}\,\mathbf z\right)+\alpha I,
	\end{align}
where 
	$$I= \displaystyle\int_\Gamma
	\left(y\cdot \nabla z-z\cdot \nabla y\right)\cdot \tau  \left(z\cdot g\right)\,dS.$$
Extending the exterior normal $n$ (defined a priori only on the boundary $\Gamma$) inside $\Omega$ by a vector field still denoted by $n$, using the Green formula and standard calculation, we can prove that for every $w\in H^2(\Omega)$ we have
	\begin{align}\label{green_convective}
	&\displaystyle\int_{\Gamma}   \left(z\cdot \nabla y-y\cdot \nabla z\right)\cdot w\,dS\notag\\
	&
	=\displaystyle\left(z\cdot \nabla y-y\cdot \nabla z,
	w\,\mathrm{div} \, n\right)+
	\int_\Omega n\cdot \nabla\left( 
	\left(z\cdot \nabla y-y\cdot \nabla z\right)\cdot w\right) \,dx\notag\\
	&=\displaystyle\left(z\cdot \nabla y-y\cdot \nabla z,
	w\,\mathrm{div} \, n\right)+b\left(n,w,z\cdot \nabla y-y\cdot \nabla z\right)+\left(\nabla
	\left(z\cdot \nabla y-y\cdot \nabla z\right)n,w\right)
	\notag\\
&=\displaystyle\left(z\cdot \nabla y-y\cdot \nabla z,
	w\,\mathrm{div} \, n\right)+b\left(n,w,z\cdot \nabla y-y\cdot \nabla z\right)+\left(\nabla
	\left(z\cdot \nabla y-y\cdot \nabla z\right),w\otimes n\right)
	\notag\\
&=\displaystyle\left(z\cdot \nabla y-y\cdot \nabla z,
	w\,\mathrm{div} \, n\right)+b\left(n,w,z\cdot \nabla y-y\cdot \nabla z\right)\notag\\
	& \ +b\left(z,\nabla y,w\otimes n\right)-b\left(y,\nabla z,w\otimes n\right)+\left(\nabla z\nabla y-\nabla y\nabla z,w\otimes n\right)\notag\\
&=\displaystyle b\left(z,y,w \,\mathrm{div} \, n\right)-
	b\left(y,z,w \,\mathrm{div} \, n\right)+b\left(n,w,z\cdot \nabla y-y\cdot \nabla z\right)\nonumber\\
	&\ -b\left(z,w\otimes n,\nabla y\right)+
	b\left(y,w\otimes n,\nabla z\right)+\left(\nabla z\nabla y-\nabla y\nabla z,w\otimes n\right),
	\end{align}
where 
	$$w\otimes n=\left(w_in_j\right)_{ij} \qquad \mbox{and} \qquad 
	b(\phi,\eta,\zeta)=\displaystyle \sum_{i,j,k}\int_\Omega \phi_k\tfrac{\partial \eta_{ij}}{\partial x_k}\zeta_{ij}\,dx$$
 Combining (\ref{trilin_0}) and (\ref{green_convective}), we deduce that
	\begin{align}\label{trilin_1}
	\left( \mathbf{curl}\,\sigma(\mathbf z)\times  \mathbf y,\mathbf z\right)
	&=b(z,y,z)+\alpha b\left(\mathbf z, \mathbf{curl}\, \mathbf y,\mathbf{curl}\, \mathbf z\right)+2\displaystyle\alpha\sum_{k=1}^2\left(\nabla \mathbf z_k\times \nabla \mathbf y_k,\mathbf{curl}\,\mathbf z\right)\notag\\
	&-\displaystyle \alpha b\left(z,y,\left(z\cdot g\right)\tau \,\mathrm{div} \, n\right)+\alpha\,
	b\left(y,z,\left(z\cdot g\right)\tau \,\mathrm{div} \, n\right)\notag\\
	&-\alpha b\left(n,\left(z\cdot g\right)\tau,z\cdot \nabla y-y\cdot \nabla z\right)
	- \alpha b\left(y,\left(z\cdot g\right)\tau\otimes n,\nabla z\right)\notag\\
	&+\alpha b\left(z,\left(z\cdot g\right)\tau\otimes n,\nabla y\right)
	+\alpha\left(\nabla z\nabla y-\nabla y\nabla z,\left(z\cdot g\right)\tau\otimes n\right).
	\end{align}
Standard arguments together with the Sobolev and the Korn inequalities give
	\begin{align}\label{trilin_2}&\left|b(z,y,z)\right|+\alpha\left|b\left(\mathbf z, \mathbf{curl}\, \mathbf y,\mathbf{curl}\, \mathbf z\right)\right|+\displaystyle 2\alpha\sum_{k=1}^2\left|\left(\nabla \mathbf z_k\times \nabla \mathbf y_k,\mathbf{curl}\,\mathbf z\right)\right|\nonumber\\
	&\leq \|z\|_4^2 \left\|\nabla y\right\|_2+
	\alpha\left(\|z\|_4 \left\|\nabla \mathbf{curl}\,\mathbf y\right\|_{4} \left\|\mathbf{curl}\, \mathbf z\right\|_2+
	2\sum_{k=1}^2\left\|\nabla \mathbf z_k\right\|_2
	\left\|\nabla \mathbf y_k\right\|_\infty\left\|\mathbf{curl}\,
	\mathbf z\right\|_2\right) \nonumber\\
	&\leq \left( \kappa_1\left\|\nabla y\right\|_2+c\,\alpha 
	\left\|y\right\|_{H^3}\right) \|\nabla z\|_2^2
	\end{align}	
and
	\begin{align}	\label{trilin_3}&\left|b\left(z,y, \left(z\cdot g\right)\tau \,\mathrm{div} \, n\right)-b\left(y,z, \left(z\cdot g\right)\tau \,\mathrm{div} \, n\right)\right|\nonumber\\
	&\leq \left(\left\|z\right\|_4\left\|\nabla y\right\|_2+
	\left\|y\right\|_4\left\|\nabla z\right\|_2
	\right)\left\| \left(z\cdot g\right)\tau \,\mathrm{div} \, n\right\|_4 \notag\\
	&\leq \left(\left\|z\right\|_4\left\|\nabla y\right\|_2+
	\left\|y\right\|_4\left\|\nabla z\right\|_2
	\right)\left\|z\right\|_4\left\|g\right\|_\infty
	\left\|\tau\right\|_\infty \left\|\mathrm{div} \, n\right\|_\infty\notag\\
	&\leq c\left\|n\right\|_\infty^2\left\|\nabla n\right\|_\infty^2\left\|\nabla y\right\|_2\left\|\nabla z\right\|_2^2\notag\\
	&\leq cc_1(n)^4\left\|y\right\|_{H^3}\left\|\nabla z\right\|_2^2
	 \end{align}
with $c_k(n)=\|n\|_{C^k(\bar\Omega)}$ and $c$ only depending on $\Omega$. Similarly, we have
	\begin{align}	\label{trilin_4}
	&\left| b\left(n,\left(z\cdot g\right)\tau,z\cdot \nabla y-y\cdot \nabla z\right)+b\left(y,\left(z\cdot g\right)\tau\otimes n,\nabla z\right)-b\left(z,\left(z\cdot g\right)\tau\otimes n,\nabla y\right)
\right|\notag\\
&\leq \left\|n\right\|_\infty
	\left(\left\|z\right\|_4\left\|\nabla y\right\|_4+
	\left\|y\right\|_\infty\left\|\nabla z\right\|_2\right)
	\left\|\nabla\left(\left(z\cdot g\right)\tau\right)\right\|_2\notag\\
&
	+\left(\left\|y\right\|_\infty\left\|\nabla z\right\|_2+
	\left\|z\right\|_4\left\|\nabla y\right\|_4 \right)\left\|\nabla\left( \left(z\cdot g\right)\tau\otimes n\right)\right\|_2
\notag\\
&\leq c\left\|y\right\|_{H^3}\left\|\nabla z\right\|_2 \left(\left\|n\right\|_\infty\left\|\nabla\left(\left(z\cdot g\right)\tau\right)\right\|_2+\left\|\nabla\left( \left(z\cdot g\right)\tau\otimes n\right)\right\|_2\right)
\notag\\
&\leq c\left\|y\right\|_{H^3}\left\|\nabla z\right\|_2 
	\left\|\nabla\left(z\cdot g\right)\right\|_2
	\left(\|n\|_\infty\|\tau\|_\infty+\left\|\tau\otimes n\right\|_\infty\right)\notag\\
&+c\left\|y\right\|_{H^3}\left\|\nabla z\right\|_2 
\left\|z\cdot g\right\|_2\left(\|n\|_\infty
	\left\|\nabla\tau\right\|_\infty+
	\left\|\nabla\left(\tau\otimes n\right)\right\|_\infty\right)
	\notag\\
	&\leq c\left\|y\right\|_{H^3}
	\left\|n\right\|_\infty^2\left( \left\|n\right\|_\infty\left\|\nabla n\right\|_\infty+\left\|\nabla n\right\|_\infty^2
	+\left\|n\right\|_\infty\left\|\nabla^{(2)} n\right\|_\infty \right)\left\|\nabla z\right\|_2^2
	\notag\\
	&\leq c c_2(n)^4 \left\|y\right\|_{H^3}\left\|\nabla z\right\|_2^2
		\end{align}
and 
\begin{align}\label{trilin5}\left|\left(\nabla z\nabla y-\nabla y\nabla z,\left(z\cdot g\right)\tau\otimes n\right)\right|&\leq 2 \left\|\nabla z\right\|_2  \left\|\nabla y\right\|_\infty  \left\|\left(z\cdot g\right)\tau\otimes n\right\|_2
	\notag\\
	&\leq c\left\|n\right\|_\infty^3
	 \left\|\nabla n\right\|_\infty 
	\left\|\nabla y\right\|_\infty \left\|\nabla z\right\|_2
	\left\|z\right\|_2\notag\\
	&\leq cc_1(n)^4
	\left\|y\right\|_{H^3} \left\|\nabla z\right\|_2^2.
\end{align}
The claimed result follows then by combining (\ref{trilin_1})-(\ref{trilin_4}).$\hfill\Box$\vspace{2mm}\\

\end{document}